\documentclass[10pt,reqno]{amsart}
\usepackage{amsmath,amssymb,latexsym,esint,cite,mathrsfs}
\usepackage{verbatim,cite,relsize,color,soul}
\usepackage[latin1]{inputenc}
\usepackage{microtype}
\usepackage{color,enumitem,graphicx}
\usepackage[colorlinks=true,urlcolor=blue, citecolor=red,linkcolor=blue,
linktocpage,pdfpagelabels, bookmarksnumbered,bookmarksopen]{hyperref}
\usepackage[hyperpageref]{backref}
\usepackage[english]{babel}
\usepackage{tikz}
\usepackage[font=small,skip=5pt]{caption}
\usepackage{geometry}
\numberwithin{equation}{section}
\newtheorem{theorem}{Theorem}[section]
\newtheorem{lemma}[theorem]{Lemma}

\newtheorem{proposition}[theorem]{Proposition}

\theoremstyle{definition}
\newtheorem{definition}[theorem]{Definition}
\newtheorem{remark}[theorem]{Remark}

\newcommand{\N}{\mathbb N}
\newcommand{\Z}{\mathbb Z}

\newcommand{\R}{\mathbb R}
\newcommand{\C}{\mathbb C}
\newcommand{\T}{\mathbb T}

\newcommand{\Pc}{\mathcal P}
\newcommand{\Zc}{\mathcal Z}
\newcommand{\Nc}{\mathcal N}

\newcommand{\vareps}{\varepsilon}

\DeclareMathOperator*{\supp}{supp}
\DeclareMathOperator*{\betc}{{\beta_c}}
\DeclareMathOperator*{\gamc}{{\gamma_c}}

\DeclareMathOperator*{\rc}{{rc}}

\newcommand{\scal}[1]{\left\langle #1 \right\rangle}

\title[Random data cubic 4NLS]
{Random data theory for the cubic fourth-order nonlinear Schr\"odinger equation}
\author[V. D. Dinh]{Van Duong Dinh}
\address[V. D. Dinh]{Laboratoire Paul Painlev\'e UMR 8524, Universit\'e de Lille CNRS, 59655 Villeneuve d'Ascq Cedex, France
and 
Department of Mathematics, HCMC University of Education, 280 An Duong Vuong, Ho Chi Minh, Vietnam}
\email{contact@duongdinh.com}


\begin{document}
	
	\begin{abstract}
	We consider the cubic nonlinear fourth-order Schr\"odinger equation
	\[
	i \partial_t u - \Delta^2 u + \mu \Delta u = \pm |u|^2 u, \quad \mu \geq 0
	\]
	on $\R^N, N\geq 5$ with random initial data. We prove almost sure local well-posedness below the scaling critical regularity. We also prove probabilistic small data global well-posedness and scattering. Finally, we prove the global well-posedness and scattering with a large probability for initial data randomized on dilated cubes.
	\end{abstract}

	\maketitle

	\section{Introduction}
	\label{S1}
	\setcounter{equation}{0}
	\subsection{Introduction}
	We consider the Cauchy problem for the cubic fourth-order nonlinear Schr\"odinger equation 
	\begin{equation} \label{4NLS}
		\left\{ 
		\begin{array}{rcl}
			i\partial_t u - \Delta^2 u + \mu \Delta u &=& \pm |u|^2 u, \quad (t,x) \in \R \times \R^N, \\
			\left.u\right|_{t=0}&=& u_0,
		\end{array}
		\right.
	\end{equation}
	where $u: \mathbb{R} \times \mathbb{R}^N \rightarrow \mathbb{C}$, $u_0: \mathbb{R}^N \rightarrow \mathbb{C}$, and $\mu \geq 0$. The plus (resp. minus) sign in front of the nonlinearity corresponds to the defocusing (resp. focusing) case. The fourth-order Schr\"odinger equation has been introduced by Karpman \cite{Karpman} and Karpman-Shagalov \cite{KS} to take into account the role of small fourth-order dispersion terms in the propagation of intense laser beams in a bulk medium with Kerr nonlinearity. 
	
	It is well-known that the equation \eqref{4NLS} with $\mu=0$ enjoys the scaling invariance
	\begin{align} \label{scaling}
	u_\lambda(t,x):= \lambda^2 u(\lambda^4 t, \lambda x), \quad \lambda>0.
	\end{align}
	A direct computation shows
	\[
	\|u_\lambda(0)\|_{\dot{H}^\gamma} = \lambda^{\gamma+2-\frac{N}{2}} \|u_0\|_{\dot{H}^\gamma}.
	\]
	We thus define the critical exponent
	\begin{align} \label{defi-gamc}
	\gamc:= \frac{N-4}{2}.
	\end{align}
	For initial data $u_0 \in H^\gamma(\R^N)$, we say that the Cauchy problem \eqref{4NLS} is subcritical, critical or supercritical if $\gamma >\gamc$, $\gamma=\gamc$ or $\gamma <\gamc$ respectively.

	\subsection{Known results}~
	In the last decade, the fourth-order Schr\"odinger equation has been attracted a lot of interest in mathematics, numerics and physics. Let us recall some known results related to \eqref{4NLS} in both deterministic and probabilistic settings. 
	
	{\bf (1) Deterministic setting.}
	Artzi-Koch-Saut \cite{BaKS} established sharp dispersive estimates for the fourth-order Schr\"odinger operator.	Thanks to these dispersive estimates, Pausader \cite{Pausader-DPDE} showed the local well-posedness for fourth-order nonlinear Schr\"odinger equations in the (sub)critical cases. Pausader \cite{Pausader-DPDE, Pausader-ener-focu, Pausader-cubic} and Miao-Xu-Zhao \cite{MXZ-1, MXZ-2} investigated the asymptotic behavior of global $H^2$-solutions in the energy-critical case. In the mass and energy intercritical case, the energy scattering for the defocusing problem was shown by Pausader \cite{Pausader-DPDE} in dimensions $N\geq 5$ and Pausader-Xia \cite{PX} in dimensions $1\leq N\leq 4$. The energy scattering for the focusing problem was studied by Guo \cite{Guo} and the author \cite{Dinh-focus}. In the mass-critical case, the asymptotic behavior of global $L^2$-solutions for was proved by Pausader-Shao \cite{PS} in dimensions $N\geq 5$. The asymptotic behavior of global solutions below the energy space was studied by Miao-Wu-Zhang \cite{MWZ} and the author \cite{Dinh-NA}. In \cite{BL}, Boulenger-Lenzmann established the existence of finite time blow-up $H^2$-solutions for the focusing problem. Dynamical properties such as mass-concentration and limiting profile of blow-up $H^2$-solutions were studied by Zhu-Yang-Zhang \cite{ZYZ-DPDE} and the author \cite{Dinh-JDDE}.
	
	{\bf (2) Probabilistic setting.}
	In the supercritical case, it was shown in \cite{Pausader-DPDE, Dinh-fourth} that \eqref{4NLS} is ill-posed in the sense that the solution map fails to be continuous at 0. Recently, the probabilistic techniques have been exploited to show almost sure local well-posedness for nonlinear dispersive equations below the critical regularity threshold. This approach was initiated by Bourgain \cite{Bourgain-CMP}. More precisely, he considered random initial data of the form
	\[
	f^\omega(x) =\sum_{n \in \Z^2} \frac{g_n(\omega)}{\sqrt{1+|n|^2}} e^{in\cdot x},
	\]
	where $(g_n)_{n\in \Z^2}$ is a sequence of independent standard complex-valued Gaussian random variables. By combining deterministic PDE techniques and probabilistic arguments, he showed that the (Wick ordered) cubic NLS on $\T^2$ with initial data $f^\omega$ is almost sure well-posed. Later, Burq-Tzvetkov \cite{BT} considered a more general class of random initial data on compact Riemannian manifolds $M$ of the form
	\[
	f^\omega(x) =\sum_{n=1}^\infty g_n(\omega) c_n e_n(x), \quad c_n = \scal{u,e_n}_{L^2(M)} = \int_M u(x) \overline{e}_n(x) d\text{vol}(x),
	\]
	where $(e_n)_{n\geq 1}$ is an orthonormal basis of $L^2(M)$ consisting of eigenfunctions of the Laplace-Beltrami operator, and $(g_n)_{n\geq 1}$ is a sequence of independent mean-zero random variables with uniform bound on the fourth moments. They proved the almost sure local well-posedness for the cubic nonlinear wave equation on the three-dimensional compact manifolds. After, Burq-Thomann-Tzvetkov \cite{BTT} and Deng \cite{Deng} showed the almost sure well-posedness for nonlinear Schr\"odinger equation with harmonic potential in dimensions 1 and 2 respectively. Recently, B\'enyi-Oh-Pocovnicu \cite{BOP-TAMS} and L\"uhrmann-Mendelson \cite{LM} independently introduced randomizations on $\R^N$. As consequences, the almost sure well-posedness for cubic nonlinear Schr\"odinger equation was shown in \cite{BOP-TAMS}, and the almost sure well-posedness for energy subcritical nonlinear wave equations was established in \cite{LM}. There are several works on random Cauchy theory followed these results (see e.g. \cite{HO, BOP, OOP, BOP-high, KMV, DLM}). 
	
	Concerning the random data Cauchy problem for fourth-order nonlinear Schr\"odinger equations, we mention recent works of Hirayama-Okamoto \cite{HO-fourth}, Chen-Zhang \cite{CZ} and Zhang-Xu \cite{ZX}. Motivated by aforementioned results, in this paper, we study the random data Cauchy problem for \eqref{4NLS}. Before stating our results, let us recall the definition of Wiener randomization on $\R^N$ due to \cite{BOP-TAMS}. Let $\psi \in C^\infty_0(\R^N)$ be such that $0\leq \psi \leq 1$, $\supp(\psi) \subset [-1,1]^N$ and 
	\begin{align} \label{defi-psi}
	\sum_{n\in \Z^N} \psi(\xi-n) =1, \quad \forall \xi \in \R^N.
	\end{align}
	Given a function $f$ on $\R^N$, we have
	\[
	f(x) = \sum_{n\in \Z^N} \psi(D-n) f(x),
	\]
	where 
	\[
	\psi(D-n)f(x) = (2\pi)^{-N} \int e^{ix\cdot \xi} \psi(\xi-n) \hat{f}(\xi) d\xi.
	\]
	The Wiener randomization of $f$ on $\R^N$ is defined by
	\begin{align} \label{defi-random}
	f^\omega(x)=\sum_{n\in \Z^N} g_n(\omega) \psi(D-n) f(x),
	\end{align}
	where $(g_n)_{n\in \Z^N}$ is a sequence of independent mean-zero complex-valued random variables on a probability space $(\Omega, \mathscr{F}, \Pc)$, where the real and imaginary parts of $g_n$ are independent and endowed with probability distributions $\mu_n^{(1)}$ and $\mu_n^{(2)}$.
	
	In the sequel, we make the following assumption: there exists $c>0$ such that
	\begin{align} \label{cond-distri}
	\left| \int_{\R} e^{\delta x} d \mu_n^{(j)}(x) \right| \leq e^{c\delta^2}
	\end{align}
	for all $\delta \in \R$, $n \in \Z^N$ and $j=1,2$. Note that \eqref{cond-distri} is satisfied by standard complex-valued Gaussian random variables, standard Bernoulli random variables and any random variables with compactly supported distributions.
	
	\subsection{Main results} 
	Denote $U_\mu(t)=e^{-it(\Delta^2-\mu \Delta)}$ the Schr\"odinger operator associated to \eqref{4NLS} and define
	\begin{align} \label{defi-gam-N}
	\gamma_N:= \max \left\{ \frac{(N-1)(N-4)}{2(N+5)}, \frac{N-4}{4}\right\}.
	\end{align}
	
	Our first result is the following almost sure local well-posedness. 
	\begin{theorem} [Almost sure local well-posedness] \label{theo-almo-sure-lwp}
		Let $N\geq 5$, $\mu \geq 0$ and $\gamma \in (\gamma_N, \gamc)$. Let $f\in H^\gamma(\R^N)$ and $f^\omega$ be the Wiener randomization defined in \eqref{defi-random} satisfying \eqref{cond-distri}. Then the equation \eqref{4NLS} is almost surely locally well-posed with respect to the randomization data $f^\omega$. More precisely, there exist $C,c, \theta >0$ such that for each $0<T\ll 1$, there exists a set $\Omega_T \subset \Omega$ with the following properties:
		\begin{itemize}
			\item $\Pc(\Omega \backslash \Omega_T) \leq C \exp \left(-c T^{-\theta} \|f\|^{-2}_{H^\gamma(\R^N)}\right)$.
			\item For each $\omega \in \Omega_T$, there exists a unique solution $u$ to \eqref{4NLS} with the initial data $f^\omega$ in the class
			\[
			U_\mu(t) f^\omega + C([0,T], H^{\gamc}(\R^N)) \subset C([0,T], H^\gamma(\R^N)).
			\]
		\end{itemize}
	\end{theorem}
	
	We will prove Theorem $\ref{theo-almo-sure-lwp}$ by considering the equation satisfied by the nonlinear part of $u$. More precisely, let 
	\begin{align} \label{z-omega}
	z(t)=z^\omega(t):= U_\mu(t) f^\omega, \quad v(t):= u(t) - U_\mu(t) f^\omega.
	\end{align}
	Then the equation \eqref{4NLS} with initial data $f^\omega$ is reduced to 
	\begin{equation} \label{4NLS-v}
	\left\{ 
	\begin{array}{rcl}
	i\partial_t v - \Delta^2 v + \mu \Delta v &=& \pm |v+z|^2 (v+z),  \\
	\left.v\right|_{t=0}&=& 0.
	\end{array}
	\right.
	\end{equation}
	We will prove the Cauchy problem \eqref{4NLS-v} is almost sure locally well-posed by viewing $z$ as a random forcing term. This is done by using variants of the Bourgain $X^{\gamma,b}$-spaces adapted to the $U^p$- and $V^p$-spaces introduced by Tataru, Koch and their collaborators \cite{HHK,HTT, KTV}. Since we are considering algebraic nonlinearity, we use the Littlewood-Paley decomposition to decompose the nonlinearity into dyadic pieces, and then carefully perform the case-by-case analysis. We refer the reader to Sections \ref{S3} and \ref{S4} for more details. 
	
	The next result is the probabilistic small data global well-posedness and scattering.
	\begin{theorem}[Probabilistic small data global well-posedness and scattering] \label{theo-proba-small-gwp}
		Let $N\geq 5$, $\mu \geq 0$ and $\gamma \in (\gamma_N, \gamc)$. Let $f\in H^\gamma(\R^N)$ and $f^\omega$ be the Wiener randomization defined in \eqref{defi-random} satisfying \eqref{cond-distri}. Then there exist $C,c>0$ such that for each $0<\vareps \ll 1$, there exists a set $\Omega_\vareps \subset \Omega$ with the following properties:
		\begin{itemize}
			\item $\Pc(\Omega \backslash \Omega_\vareps) \leq C \exp \left(-c \vareps^{-2} \|f\|^{-2}_{H^\gamma(\R^N)}\right) \rightarrow 0$ as $\vareps \rightarrow 0$.
			\item For each $\omega \in \Omega_\vareps$, there exists a unique global in time solution $u$ to \eqref{4NLS} with initial data $\vareps f^\omega$ in the class
			\[
			\vareps U_\mu(t) f^\omega + C(\R, H^{\gamc}(\R^N)) \subset C(\R, H^\gamma(\R^N)).
			\]
			\item For each $\omega \in \Omega_\vareps$, there exists $f^\omega_+ \in H^{\gamc}(\R^N)$ such that 
			\[
			\|u(t)-\vareps U_\mu(t) f^\omega - U_\mu(t) f^\omega_+\|_{H^{\gamc}(\R^N)} \rightarrow 0
			\]
			as $t\rightarrow \infty$. A similar statement holds for $t\rightarrow -\infty$.
		\end{itemize}
	\end{theorem}
	
	\begin{remark}
		In \cite{CZ}, Chen-Zhang considered the Cauchy problem of the fourth-order nonlinear Schr\"odinger equation of the form
		\[
		i \partial_t u +\mu \Delta u +\Delta^2 u = P_m\left((\partial^\alpha_x u)_{|\alpha| \leq 2}, (\partial^\alpha_x \overline{u})_{|\alpha| \leq 2} \right), 
		\]
		where $\mu \in \R$, $P_m$ is a homogeneous polynomial of degree $m\geq 3$ containing the second order derivative. More precisely, $P_m$ is of the form
		\[
		c_1 \prod_{k=1}^m u_k + c_2 \sum_{i=1}^N \partial_i u_1 \prod_{k=2}^m u_k + c_3 \sum_{i,j=1}^N \partial_i u_1 \partial_j u_2 \prod_{k=3}^m u_k + c_4 \sum_{i,j=1}^N \partial^2_{ij} u_1 \prod_{k=2}^m u_k,
		\]
		where $u_k=u$ or $\overline{u}$ for $k=1,\cdots,m$. They established almost sure local well-posedness and probabilistic small data global existence and scattering for random initial data in $H^\gamma$ with $\gamma \in (\beta_{N,m}, \beta_{\text{c},m}]$, where
		\begin{align*}
		\beta_{\text{c},m}:= \frac{N}{2} -\frac{2}{m-1}, \quad \beta_{N,m}:= \left\{
		\renewcommand*{\arraystretch}{1.2}
		\begin{array}{lcl}
		\beta_{\text{c},m} -\frac{1}{2} +\frac{m-2}{3m-7} &\text{if}& N=2, m\geq 4, \\
		\beta_{\text{c},m} -\frac{1}{2} +\frac{5-m}{2(N-1)(m-1)} &\text{if}& N\geq 3, 3\leq m <5, \\
		\beta_{\text{c},m} - \frac{1}{2} &\text{if}& N\geq 3, m\geq 5.
		\end{array}
		\right.
		\end{align*}
		In particular, when $m=3$, we have
		\[
		\betc:= \beta_{N,3} = \frac{N-2}{2}, \quad \beta_N:= \beta_{\text{c},3} = \frac{(N-2)^2}{2(N-1)}.
		\]
		It is easy to see that $\beta_N > \gamc$. Thus, the result in \cite{CZ} does not apply to show almost sure well-posedness for \eqref{4NLS} below the critical regularity threshold.
	\end{remark}
	
	Finally, we have the almost sure global well-posedness and scattering with a large probability. This is done by considering the randomization based on a partition of the frequency space by dilated cubes. More precisely, given $\lambda>0$, we define 
	\[
	\psi_\lambda(\xi) := \psi(\lambda \xi),
	\] 
	where $\psi$ is as in \eqref{defi-psi}. We can write a function $f$ on $\R^N$ as
	\[
	f(x) = \sum_{n\in \Z^N} \psi_\lambda(D-\lambda^{-1} n) f(x).
	\]
	We now introduce the randomization $f^{\omega,\lambda}$ of $f$ on dilated cubes of scale $\lambda$ by
	\begin{align} \label{defi-rando-lambda}
	f^{\omega,\lambda}(x):= \sum_{n\in \Z^N} g_n(\omega) \psi_\lambda(D-\lambda^{-1} n) f(x),
	\end{align}
	where $(g_n)_{n\in \Z^N}$ is a sequence of independent mean-zero complex-valued random variables satisfying \eqref{cond-distri}. We have the following global well-posedness and scattering with a large probability.
	
	\begin{theorem} [Large probability global well-posedness and scattering] \label{theo-gwp-large}
		Let $N\geq 5$, $\mu = 0$ and $\gamma \in (\gamma_N, \gamc)$. Let $f\in H^\gamma(\R^N)$ and $f^{\omega,\lambda}$ be the Wiener randomization on dilated cubes of scale $\lambda \gg 1$ defined in \eqref{defi-rando-lambda}. Then the equation \eqref{4NLS} is globally well-posed with a large probability. More precisely, for each $0<\vareps \ll 1$, there exists a large dilation scale $\lambda_0 = \lambda_0(\vareps, \|f\|_{H^\gamma}) >0$ such that for each $\lambda >\lambda_0$, there exists a set $\Omega_\lambda \subset \Omega$ with the following properties:
		\begin{itemize}
			\item $\Pc(\Omega \backslash \Omega_\lambda) <\vareps$.
			\item For each $\omega \in \Omega_\lambda$, there exists a unique global-in-time solution to \eqref{4NLS} with initial data $f^{\omega,\lambda}$ in the class
			\[
			U_0(t) f^{\omega,\lambda} + C(\R, H^{\gamc}(\R^N)) \subset C(\R, H^\gamma(\R^N)).
			\]
			\item For each $\omega \in \Omega_\lambda$, there exists $f^\omega_+ \in H^{\gamc}(\R^N)$ such that
			\[
			\|u(t) - U_0(t) f^{\omega,\lambda} - U_0(t) f^\omega_+\|_{H^{\gamc}(\R^N)} \rightarrow 0
			\] 
			as $t\rightarrow \infty$. A similar statement holds for $t\rightarrow -\infty$. 
		\end{itemize}
	\end{theorem}

	This paper is organized as follows. In Section $\ref{S2}$, we give some preliminaries needed in the sequel including some basic properites of the Wiener randomization, probabilistic Strichartz estimates and function spaces. In Section $\ref{S3}$, we prove probabilistic nonlinear estimates which are key ingredients of the proof. Finally, in Section $\ref{S4}$, we prove the almost sure well-posedness given in Theorems $\ref{theo-almo-sure-lwp}$, $\ref{theo-proba-small-gwp}$ and $\ref{theo-gwp-large}$. 
	
	\section{Preliminaries}
	\label{S2}
	\setcounter{equation}{0}
	
	\subsection{Notations}
	Let $1 \leq r \leq \infty$ and $\gamma \in \R$. We denote the Lebesgue space and Sobolev space by $L^r(\R^N)$ and $H^\gamma(\R^N)$ respectively. The notation $A \lesssim B$ means that there exists a constant $C>0$ such that $A \leq C B$. Similarly, $A \gtrsim B$ means $A \geq cB$ for some constant $c>0$. We also use $A \sim B$ if $A \lesssim B$ and $A \gtrsim B$. Let $I\subset \R$ and $1 \leq q,r \leq \infty$. We define the mixed norm 
	\[
	\|u\|_{L^q_tL^r_x(I\times \R^N)} := \left( \int_I \left( \int_{\R^N} |u(t,x)|^r dx \right)^{\frac{q}{r}} dt \right)^{\frac{1}{q}}
	\]
	with a usual modification when either $q$ or $r$ are infinity. When $q=r$, we use the notation $L^q_{t,x}(I \times \R^N)$ instead of $L^q_tL^q_x(I\times \R^N)$. Let $\varphi$ be a smooth real-valued radial function on $\R^N$ satisfying
	\[
	\varphi(\xi) = \left\{
	\renewcommand*{\arraystretch}{1.3}
	\begin{array}{lcl}
	1 &\text{if}& |\xi| \leq 1, \\
	0 &\text{if} & |\xi| \geq 2.
	\end{array}
	\right.
	\]
	We define the Littlewood-Paley operators
	\[
	\widehat{P_1 f} (\xi) = \varphi(\xi) \hat{f}(\xi)
	\]
	and for dyadic numbers $M= 2^m, m\geq 1$, 
	\[
	\widehat{P_M f}(\xi) = \left( \varphi(M^{-1}\xi) - \varphi(2M^{-1}\xi)\right) \hat{f}(\xi).
	\]
	We also define
	\[
	P_{\leq M_1} =\sum_{1\leq M \leq M_1} P_M, \quad P_{\geq M_1} = \sum_{M\geq M_1} P_M.
	\]
	
	\subsection{Wiener randomization and Probabilistic Strichartz estimates}
	In this subsection, we first recall some basic properties of the Wiener randomization and probabilistic Strichartz estimates related to the fourth-order Schr\"odinger equation.
	
	The first property of the Wiener randomization is that it preserves the differentiability in the sense: if $f \in H^\gamma(\R^N)$, then $f^\omega \in H^\gamma(\R^N)$ almost surely. 
	\begin{lemma}[{\cite[Lemma 3]{BOP}}] \label{lem-differ}
		Let $f\in H^\gamma(\R^N)$ and $f^\omega$ be the Wiener randomization defined in \eqref{defi-random} satisfying \eqref{cond-distri}. Then it holds that
		\[
		\Pc\left(\|f^\omega\|_{H^\gamma(\R^N)} >\lambda\right) \leq C \exp \left(-c \lambda^2 \|f\|^{-2}_{H^\gamma(\R^N)} \right)
		\]
		for all $\lambda>0$. In particular, $f^\omega \in H^\gamma(\R^N)$ almost surely.
	\end{lemma}
	
	\begin{remark}
		It was shown in \cite[Appendix]{BT} that the Wiener randomization does not gain differentiability in the sense: if $f \in H^\gamma(\R^N) \backslash H^{\gamma+\vareps}(\R^N)$, then $f^\omega \in H^\gamma(\R^N) \backslash H^{\gamma+\vareps}(\R^N)$ almost surely. 
	\end{remark}
	
	The second property of the Wiener randomization is that it gains the integrability in the sense: if $f \in L^2(\R^N)$, then $f^\omega \in L^p(\R^N)$ for all $2 \leq p <\infty$ almost surely. 
	
	\begin{lemma} [{\cite[Lemma 4]{BOP}}] \label{lem-interga}
		Let $f \in L^2(\R^N)$ and $f^\omega$ be the Wiener randomization defined in \eqref{defi-random} satisfying \eqref{cond-distri}. Then it holds that
		\[
		\Pc \left(\|f^\omega\|_{L^p(\R^N)} >\lambda \right) \leq C \exp \left(-c\lambda^2 \|f\|^{-2}_{L^2(\R^N)}\right)
		\]
		for all $p \in [2,\infty)$ and all $\lambda>0$. In particular, $f^\omega \in L^p(\R^N), p\in [2,\infty)$ almost surely.
	\end{lemma}
	
	\begin{remark}
		Comparing to the Sobolev embedding $H^{\frac{N}{2}}(\R^N) \hookrightarrow L^p(\R^N)$ for all $p \in [2,\infty)$, the Wiener randomization makes a gain of $\frac{N}{2}$ derivatives.
	\end{remark}
	
	The Wiener randomization also allows us establish some improvements of Strichartz estimates which are essential tools to the almost sure well-posedness of \eqref{4NLS}. Let us recall Strichartz estimates for \eqref{4NLS} on $\R^N$. A pair $(q,r)$ is Biharmonic admissible, or $(q,r) \in B$ for short, if 
	\[
	\frac{4}{q}+\frac{N}{r}=\frac{N}{2}, \quad \left\{
	\renewcommand*{\arraystretch}{1.2}
	\begin{array}{lcl}
	r \in \left[2,\frac{2N}{N-4}\right] &\text{if}& N\geq 5, \\
	r \in [2,\infty) &\text{if}& N=4, \\
	r\in [2,\infty] &\text{if} & N \leq 3.
	\end{array}
	\right.
	\]
	
	Let $\mu \in \R$. We denote $U_\mu(t):= e^{-it(\Delta^2-\mu\Delta)}$ the propagator for the free fourth-order Schr\"odinger equation
	\[
	i\partial_t -\Delta^2 u +\mu \Delta u=0.
	\]
	We have the following dispersive estimates due to Ben Artzi-Koch-Saut \cite{BaKS}.
	\begin{lemma}[{\cite[Theorem 1]{BaKS}}] \label{lem-dis-est}
		Let $N\geq 1$ and $\mu \in \R$. It holds that
		\[
		\|U_\mu(t) f\|_{L^\infty(\R^N)} \lesssim |t|^{-\frac{N}{4}} \|f\|_{L^1(\R^N)}
		\]
		for all $t \ne 0$, and if $\mu <0$, it requires $|t| \leq 1$.
	\end{lemma}

	Using this dispersive estimate and the abstract theory of Keel-Tao \cite{KT}, we have the following Strichartz estimates for \eqref{4NLS}.
	\begin{lemma}[Strichartz estimates {\cite[Proposition 3.1]{Pausader-DPDE}}] \label{lem-str-est}
		Let $N\geq 1$, $\mu \in \R$ and $I \subset \R$ be an interval. It holds that
		\[
		\|U_\mu(t) f\|_{L^q_tL^r_x(I \times \R^N)} \lesssim \|f\|_{L^2(\R^N)}
		\]
		for all Biharmonic admissible pairs $(q,r)$, and if $\mu<0$, it requires $|I| <\infty$.
	\end{lemma}
	
	\begin{remark} 
		We have from Sobolev embedding and Strichartz estimates that
		\begin{align} \label{est-Lp}
		\|U_\mu(t) f\|_{L^p_{t,x}(I \times \R^N)} \lesssim \||\nabla|^{\frac{N}{2} - \frac{N+4}{p}} f\|_{L^2(\R^N)}
		\end{align}
		for $p\geq \frac{2(N+4)}{N}$. Note that the derivative loss in \eqref{est-Lp} depends only on the size of the frequency support and not its location. Namely, if $\hat{f}$ is supported on a cube $Q$ of side length $M$, then
		\[
		\|U_\mu(t)f\|_{L^p_{t,x}(I\times \R^N)} \lesssim M^{\frac{N}{2}-\frac{N+4}{p}} \|f\|_{L^2(\R^N)}
		\]
		which follows from Bernstein's inequalities.
	\end{remark}
	
	We have the following improvements of Strichartz estimates under the Wiener randomization.
	\begin{lemma} [Local-in-time probabilistic Strichartz estimates] \label{lem-local-proba-str-est}
		Let $N\geq 1$ and $\mu \in \R$. Let $f \in L^2(\R^N)$ and $f^\omega$ be the Wiener randomization defined in \eqref{defi-random} satisfying \eqref{cond-distri}. Then for $2\leq q,r <\infty$, there exist $C,c>0$ such that 
		\[
		\Pc \left( \|U_\mu(t) f^\omega\|_{L^q_t L^r_x([0,T]\times \R^N)} >\lambda \right) \leq C \exp \left(-c \lambda^2 T^{-\frac{2}{q}} \|f\|^{-2}_{L^2(\R^N)}\right)
		\]
		for all $T>0$ and all $\lambda>0$. 
	\end{lemma}
	
	\begin{remark}
		Taking $\lambda=T^\theta R$, we have
		\begin{align} \label{local-proba-str-est}
		\|U_\mu(t) f^\omega\|_{L^q_t L^r_x([0,T]\times \R^N)} \lesssim T^\theta R
		\end{align}
		outside a set of probability at most 
		\[
		C \exp \left(- cR^2 T^{-2\left(\frac{1}{q}-\theta\right)} \|f\|^{-2}_{L^2(\R^N)} \right)
		\]
		for all $T, \theta, R>0$. Note that for $R>0$ fixed, this probability can be made arbitrarily small by letting $T\rightarrow 0$ as long as $\theta <\frac{1}{q}$.
	\end{remark}
	
	\begin{lemma} [Global-in-time probabilistic Strichartz estimates] \label{lem-glo-proba-str-est}
		Let $N\geq 1$ and $\mu \geq 0$. Let $f \in L^2(\R^N)$ and $f^\omega$ be the Wiener randomization defined in \eqref{defi-random} satisfying \eqref{cond-distri}. Let $(q,r)$ be a Biharmonic admissible pair with $q,r <\infty$. Let $\tilde{r} \geq r$. Then there exist $C,c>0$ such that
		\[
		\Pc \left(\|U_\mu(t)f^\omega\|_{L^q_tL^{\tilde{r}}_x(\R \times \R^N)} >\lambda \right) \leq C \exp \left(-c \lambda^2 \|f\|^{-2}_{L^2(\R^N)} \right)
		\]
		for all $\lambda>0$. 
	\end{lemma}
	
	The proofs of Lemmas $\ref{lem-local-proba-str-est}$ and $\ref{lem-glo-proba-str-est}$ follow the same argument as in \cite{BOP} using Strichartz estimates given in Lemma $\ref{lem-str-est}$. We thus omit the details.

	\subsection{Function spaces and their properties}
	In this subsection, we recall the definitions and basic properties of the $U^p$- and $V^p$-spaces developed by Tataru, Koch and their collaborators \cite{HHK, HTT, KTV}. These spaces have been very effective in establishing well-posedness of various dispersive PDEs in critical regularities. 
	
	Let $\Zc$ be the the collection of finite partitions $(t_k)_{k=0}^K$ of $\R$, i.e. $-\infty <t_0<\cdots < t_K \leq \infty$. If $t_K=\infty$, we use the convention $u(t_K):= 0$ for all functions $u: \R \rightarrow H^\gamma(\R^N)$.
	
	\begin{definition}
		Let $1\leq p<\infty$ and $\gamma \in \R$. 
		\begin{itemize}
			\item A $U^p$-atom is defined by a step function $a:\R \rightarrow H^\gamma(\R^N)$ of the form
			\[
			a(t)= \sum_{k=1}^K \phi_{k-1} \chi_{[t_{k-1},t_k)}(t), \quad \sum_{k=0}^{K-1} \|\phi_k\|_{H^\gamma(\R^N)}^p = 1,
			\]
			where $(t_k)_{k=0}^K \in \Zc$, $(\phi_k)_{k=0}^{K-1} \subset H^\gamma(\R^N)$ and $\chi_I$ is the characteristic function of $I$. 
			\item We define the atomic space $U^p(\R,H^\gamma(\R^N))$ to be the collection of functions $u:\R \rightarrow H^\gamma(\R^N)$ of the form
			\begin{align} \label{defi-u}
			u= \sum_{j=1}^\infty \lambda_j a_j
			\end{align}
			with the norm
			\[
			\|u\|_{U^p(\R,H^\gamma(\R^N))}:=\inf \left\{ \sum_{j=1}^\infty |\lambda_j| \ : \ \eqref{defi-u} \text{ holds}\right\},
			\]
			where $a_j$ are $U^p$-atoms and $(\lambda)_{j\geq 1} \in \ell^1(\C)$.
			\item We define the space of bounded $p$-variation $V^p(\R,H^\gamma(\R^N))$ to be the collection of functions $u:\R\rightarrow H^\gamma(\R^N)$ with the norm
			\[
			\|u\|_{V^p(\R,H^\gamma(\R^N))}:= \sup_{(t_k)_{k=0}^K\in \Zc}  \left(\sum_{k=1}^K \|u(t_k)-u(t_{k-1})\|_{H^\gamma(\R^N)}^p\right)^{\frac{1}{p}}.
			\]
			We also define $V^p_{\rc}(\R,H^\gamma(\R^N))$ to be the closed subspace of all right-continuous functions in $V^p(\R,H^\gamma(\R^N))$ such that $\lim_{t\rightarrow -\infty} u(t)=0$.
			\item We define $U^p_{\Delta} H^\gamma(\R^N)$ ( resp. $V^p_{\Delta}H^\gamma(\R^N)$) to be the space of all functions $u: \R \rightarrow H^\gamma(\R^N)$ such that the following norm is finite:
			\[
			\|u\|_{U^p_\Delta H^\gamma(\R^N)}:= \|U_\mu(-t) u\|_{U^p(\R, H^\gamma(\R^N))} \quad \left( \text{resp. } \|v\|_{V^p_\Delta H^\gamma(\R^N)}:= \|U_\mu(-t) u\|_{V^p_\Delta H^\gamma(\R^N)}\right). 
			\]
			We use $V^p_{\rc,\Delta} H^\gamma(\R^N)$ to denote the subspace of right-continuous functions in $V^p_\Delta H^\gamma(\R^N)$.
		\end{itemize}
	\end{definition}
	
	\begin{remark}
		It was shown in \cite{HHK} that the spaces $U^p(\R, H^\gamma(\R^N)), V^p(\R, H^\gamma(\R^N))$ and $V^p_{\rc}(\R, H^\gamma(\R^N))$ are Banach spaces. The closed subspace of continuous functions in $U^p(\R, H^\gamma(\R^N))$ is also a Banach space. Moreover, we have the following embeddings:
		\begin{align*} 
		U^p(\R, H^\gamma(\R^N)) \hookrightarrow V^p_{\rc}(\R, H^\gamma(\R^N)) \hookrightarrow U^q(\R, H^\gamma(\R^N)) \hookrightarrow L^\infty(\R, H^\gamma(\R^N))
		\end{align*}
		for $1\leq p<q<\infty$. Similar embeddings hold for $U^p_\Delta H^\gamma(\R^N)$ and $V^p_\Delta H^\gamma(\R^N)$.
	\end{remark}
	
	We have the following tranference principle. 
	\begin{lemma}[Tranference principle \cite{HHK}] 
		Let $N \geq 1$ and $\mu \geq 0$. Let $T$ be a $k$-linear operator. Suppose that we have
		\[
		\|T(U_\mu(t) f_1, \cdots, U_\mu(t) f_k)\|_{L^p_tL^q_x(\R\times \R^N)} \lesssim \prod_{j=1}^k \|f_j\|_{L^2(\R^N)}
		\]
		for some $1\leq p,q\leq \infty$. Then, we have
		\[
		\|T(u_1, \cdots, u_k)\|_{L^p_t L^q_x(\R \times \R^N)} \lesssim \prod_{j=1}^k \|u_j\|_{U^p_\Delta L^2(\R^N)}.
		\]
	\end{lemma}
	We also have the following interpolation inequality.
	\begin{lemma}[Interpolation lemma \cite{HHK}]
		Let $\gamma \in \R$ and $E$ be a Banach space. Suppose that $T: U^{p_1}(\R, H^\gamma(\R^N)) \times \cdots \times U^{p_k}(\R, H^\gamma(\R^N)) \rightarrow E$ is a bounded $k$-linear operator such that 
		\[
		\|T(u_1, \cdots, u_k)\|_{E} \leq C_1 \prod_{j=1}^k \|u_j\|_{U^{p_j}(\R, H^\gamma(\R^N))}
		\]
		for some $p_1, \cdots, p_k>2$. Moreover, assume that there exists \footnote{Since $U^2(\R, H^\gamma(\R^N)) \hookrightarrow U^p(\R,H^\gamma(\R^N))$, we have $C_2 \leq C_1$.} $C_2 \in (0,C_1]$ such that 
		\[
		\|T(u_1, \cdots, u_k)\|_E \leq C_2 \prod_{j=1}^k \|u_j\|_{U^2(\R, H^\gamma(\R^N))}.
		\]
		Then we have
		\[
		\|T(u_1, \cdots, u_k)\|_E \leq C_2 \left(\ln \frac{C_1}{C_2} +1\right)^k \prod_{j=1}^k \|u_j\|_{V^2(\R, H^\gamma(\R^N))}
		\]
		for $u_j \in V^2_{\rc}(\R, H^\gamma(\R^N)), j=1, \cdots, k$.
	\end{lemma}
	We refer the reader to \cite{HHK} (see also \cite{KTV}) for the proof of above results.
	
	\begin{definition}
		Let $\gamma \in \R$.
		\begin{itemize}
			\item We define $X^\gamma(\R)$ to be the space of all tempered distributions $u:\R \rightarrow H^\gamma(\R^N)$ such that the norm
			\[
			\|u\|_{X^\gamma(\R)} := \left(\sum_{M\geq 1 \atop \text{dyadic}} M^{2\gamma} \|P_M u\|^2_{U^2_\Delta L^2(\R^N)} \right)^{\frac{1}{2}}
			\]
			is finite.
			\item We define $Y^\gamma(\R)$ to be the space of all tempered distributions $u:\R \rightarrow H^\gamma(\R^N)$ such that for every dyadic number $M = 2^{m}, m\geq 0$, the map $t\mapsto P_M u(t)$ is in $V^2_{\rc,\Delta}H^\gamma(\R^N)$ and the norm
			\[
			\|u\|_{Y^\gamma(\R)}:= \left( \sum_{M\geq 1 \atop \text{dyadic}} M^{2\gamma} \|P_M u\|^2_{V^2_\Delta L^2(\R^N)} \right)^{\frac{1}{2}}
			\]
			is finite.
		\end{itemize}
	\end{definition}
	By definition, we have 
	\begin{align} \label{prop-X-gamma}
	\|U_\mu (t) f\|_{X^\gamma(\R)} \sim \|f\|_{H^\gamma(\R^N)}.
	\end{align}
	Moreover, we have the following embeddings:
	\begin{align} \label{embed}
	U^2_\Delta H^\gamma(\R^N) \hookrightarrow X^\gamma(\R) \hookrightarrow Y^\gamma(\R) \hookrightarrow V^2_\Delta H^\gamma (\R^N) \hookrightarrow U^p_\Delta H^\gamma (\R^N) \hookrightarrow L^\infty(\R, H^\gamma(\R^N))
	\end{align}
	for $p>2$. 
	
	Given an interval $I\subset \R$, we define the local-in-time versions $X^\gamma(I)$ and $Y^\gamma(I)$ of these spaces as restriction norms. For example, we define the $X^\gamma(I)$-norm by
	\[
	\|u\|_{X^\gamma(I)} := \inf \left\{ \|v\|_{X^\gamma(\R)} \ : \ v|_{I}=u\right\}.
	\]
	We have the following lemmas due to B\'enyi-Oh-Pocovnicu \cite[Appendix]{BOP-TAMS}. 
	\begin{lemma} [{\cite[Appendix]{BOP-TAMS}}] \label{lem-rest-norm}
		Let $\gamma \in \R$ and $1\leq p<\infty$. Let 
		\[
		u= \sum_{j=1}^\infty \lambda_j a_j,
		\]
		where $(\lambda_j)_{j\geq 1} \in \ell^1(\C)$ and $a_j$ are $U^p$-atoms. Given an interval $I\subset \R$. Then we can write
		\[
		u \cdot \chi_I = \sum_{j=1}^\infty \tilde{\lambda}_j \tilde{a}_j
		\]
		for some $(\tilde{\lambda}_j)_{j\geq 1} \in \ell^1(\C)$ and $\tilde{a}_j$ are $U^p$-atoms satisfying
		\[
		\sum_{j=1}^\infty |\tilde{\lambda}_j| \leq \sum_{j=1}^\infty |\lambda_j|.
		\]
		As a consequence, we have
		\[
		\|u\cdot \chi_I\|_{U^p(\R, H^\gamma(\R^N))} \leq \|u\|_{U^p(\R,H^\gamma(\R^N))}
		\]
		for any $u\in U^p(\R,H^\gamma(\R^N))$ and any $I \subset \R$.
	\end{lemma}
	
	Given an interval $I \subset \R$. We define the local-in-time $U^p$-norm in the usual manner as a restriction norm 
	\[
	\|u\|_{U^p(I,H^\gamma(\R^N))} =\inf \left\{ \|v\|_{U^p(\R, H^\gamma(\R^N))} \ : \ v|_I = u\right\}.
	\]
	Note that this infimum achieve by $v= u \cdot \chi_I$ in view of Lemma $\ref{lem-rest-norm}$. 
	
	We have the following Strichartz estimates adapted to the $X^\gamma$- and $Y^\gamma$-spaces.
	\begin{lemma} \label{lem-str-est-Y}
		Let $N\geq 1$ and $\mu \geq 0$. Let $(q,r)$ be a Biharmonic admissible pair with $q>2$ and $p\geq \frac{2(N+4)}{N}$. Then for any $0<T\leq \infty$, we have
		\begin{align}
		\|u\|_{L^q_tL^r_x([0,T)\times \R^N)} &\lesssim \|u\|_{Y^0([0,T))}, \label{est-LqLr-Y} \\
		\|u\|_{L^p_{t,x}([0,T)\times \R^N)} &\lesssim \||\nabla|^{\frac{N}{2}-\frac{N+4}{p}} u\|_{Y^0([0,T))}. \label{est-Lp-Y} 
		\end{align}
	\end{lemma}
	
	\begin{proof}
		By Strichart estimates, we have
		\[
		\|U_\mu(t) f\|_{L^q_t L^r_x ([0,T)\times \R^N)} \lesssim \|f\|_{L^2(\R^N)}.
		\]
		It follows from the tranference principle that
		\[
		\|u\|_{L^q_tL^r_x([0,T)\times \R^N)} \lesssim \|u\|_{U^q_\Delta([0,T),L^2(\R^N))} \leq \|u\|_{Y^0([0,T))},
		\]
		where we have used the embedding \eqref{embed}.
		
		Similarly, by \eqref{est-Lp} and the tranference principle, we have
		\[
		\|u\|_{L^p_{t,x}([0,T)\times \R^N)} \lesssim \||\nabla|^{\frac{N}{2}-\frac{N+4}{p}} u\|_{U^p_\Delta([0,T),L^2(\R^N))} \leq \||\nabla|^{\frac{N}{2}-\frac{N+4}{p}} u\|_{Y^0([0,T))}.
		\]
	\end{proof}
	
	\begin{remark} \label{rem-est-Y}
		\begin{itemize}
			\item The derivative loss in \eqref{est-Lp-Y} depends only on the size of the spatial frequency support and not its location. Namely, if the spatial frequency support of $\hat{u}(t,\xi)$ is contained in a cube of side length $M$ for all $t\in \R$, then 
			\begin{align*} 
			\|u\|_{L^p_{t,x}([0,T) \times \R^N)} \lesssim M^{\frac{N}{2}-\frac{N+4}{p}} \|u\|_{Y^0([0,T))}
			\end{align*}
			by Bernstein's inequalities.
			\item By \eqref{embed}, we can replace the norm $Y^0([0,T))$ in \eqref{est-LqLr-Y} and \eqref{est-Lp-Y} by $X^0([0,T))$.
		\end{itemize}
	\end{remark}
	
	We also have the following bilinear estimate related to the fourth-order Schr\"odinger equation.
	
	\begin{lemma}[Bilinear estimate] \label{lem-bili-est}
		Let $N\geq 5$ and $\mu \geq 0$. Let $M_1, M_2 \in 2^{\N}$ be such that $M_1 \leq M_2$. Then it holds that
		\begin{align} \label{bili-est}
		\|[U_\mu(t) P_{M_1} f] [U_\mu(t) P_{M_2} g] \|_{L^2_{t,x}(\R \times \R^N)} \lesssim M_1^{\frac{N-4}{2}} \left(\frac{M_1}{M_2}\right)^{\frac{3}{2}} \|f\|_{L^2(\R^N)} \|g\|_{L^2(\R^N)}.
		\end{align}
	\end{lemma}
	
	\begin{proof}
		For simplifying the notation, we denote $L^p_t L^q_x$, $L^p_{t,x}$ and $L^2_x$ instead of $L^p_tL^q_x(\R \times \R^N), L^p_{t,x}(\R \times \R^N)$ and $L^2(\R^N)$ respectively. 
		
		We first consider the case $M_1 \sim M_2$. By H\"older's inequality, Sobolev embedding and Strichartz estimates, we have
		\begin{align*}
		\|[U_\mu(t) P_{M_1}f][U_\mu(t) P_{M_2}g]\|_{L^2_{t,x}} &\lesssim \|U_\mu(t) P_{M_1} f\|_{L^4_t L^N_x} \|U_\mu(t) P_{M_2} g\|_{L^4_t L^{\frac{2N}{N-2}}_x} \\
		&\lesssim \||\nabla|^{\frac{N-4}{2}} U_\mu(t) P_{M_1} f\|_{L^4_t L^{\frac{2N}{N-2}}_x} \|U_\mu(t) P_{M_2} g\|_{L^4_t L^{\frac{2N}{N-2}}_x} \\
		&\lesssim M_1^{\frac{N-4}{2}}  \|f\|_{L^2_x} \|g\|_{L^2_x}\sim M_1^{\frac{N-4}{2}} \left(\frac{M_1}{M_2}\right)^{\frac{3}{2}} \|f\|_{L^2_x} \|g\|_{L^2_x}.
		\end{align*}
		Here we have used the fact $\left(4,\frac{2N}{N-2}\right)$ is Biharmonic admissible.
			
		We next consider the case $M_1 \ll M_2$. By duality, we have
		\[
		\text{LHS} \eqref{bili-est} = \sup_{\|G\|_{L^2_{t,x}} =1} \left| \scal{[U_\mu(t) P_{M_1} f] [U_\mu(t) P_{M_2} g], G}_{L^2_{t,x}} \right|,
		\]
		where
		\[
		\scal{F,G}_{L^2_{t,x}} = \iint_{\R \times \R^N} F(t,x)\overline{G}(t,x) dx dt = \int_{\R} \scal{F(t), G(t)}_{L^2_x} dt.
		\]
		By Parseval's identity, we have
		\[
		\text{LHS} \eqref{bili-est} = \sup_{\|G\|_{L^2_{t,x}} =1} \left| \int_{\R} \scal{ \mathcal{F}\left( [U_\mu(t)P_{M_1} f] [U_\mu(t) P_{M_2}g] \right), \hat{G}(t) }_{L^2_\xi} dt\right|,
		\]
		where
		\begin{align*}
		\mathcal{F} \left( [U_\mu(t)P_{M_1} f] [U_\mu(t) P_{M_2}g] \right) (\xi) &= \int_{\R^N} e^{-it \left( |\xi-\eta|^4- \mu |\xi-\eta|^2\right)} \widehat{P_{M_1} f} (\xi-\eta) e^{-it\left( |\eta|^4-\mu|\eta|^2\right)} \widehat{P_{M_2} g}(\eta) d\eta \\
		&=\int_{\R^N} e^{-it \left((|\xi-\eta|^4+|\eta|^4) - \mu(|\xi-\eta|^2+ |\eta|^2) \right)} \widehat{P_{M_1} f} (\xi-\eta) \widehat{P_{M_2} g}(\eta) d\eta.
		\end{align*}
		It follows that
		\begin{align*}
		\int_{\R} &\scal{ \mathcal{F}\left( [U_\mu(t)P_{M_1} f] [U_\mu(t) P_{M_2}g] \right), \hat{G}(t)}_{L^2_\xi} dt \\
		&=\int_{\R} \scal{\int_{\R^N} e^{-it \left(|\xi-\eta|^4+|\eta|^4 - \mu(|\xi-\eta|^2+ |\eta|^2) \right)} \widehat{P_{M_1} f} (\xi-\eta) \widehat{P_{M_2} g}(\eta) d\eta, \hat{G}(t) }_{L^2_\xi} dt \\
		&= \int_{\R^N} \scal{\widehat{P_{M_1}f}(\cdot-\eta) \widehat{P_{M_2}g}(\eta), \tilde{G}\left( |\cdot-\eta|^4 + |\eta|^4- \mu (|\cdot-\eta|^2+|\eta|^2), \cdot\right)}_{L^2_\xi} d\eta \\
		&= \iint_{\R^N\times \R^N}  \widehat{P_{M_1}f}(\xi-\eta) \widehat{P_{M_2}g}(\eta) \overline{\tilde{G}} \left( |\xi-\eta|^4 + |\eta|^4- \mu (|\xi-\eta|^2+|\eta|^2), \xi\right) d\xi d\eta,
		\end{align*}
		where $\tilde{G}$ is the space-time Fourier transform of $G$. Thus, the estimate \eqref{bili-est} is reduced to show
		\begin{multline} \label{bili-est-proof}
		\left|\iint_{\R^N\times \R^N} \overline{\tilde{G}} \left(|\xi|^4+|\eta|^4 - \mu(|\xi|^2+|\eta|^2), \xi+\eta\right) \widehat{P_{M_1} f}(\xi) \widehat{P_{M_2} g}(\eta) d\xi d\eta\right| \\
		\lesssim M_1^{\frac{N-4}{2}} \left(\frac{M_1}{M_2}\right)^{\frac{3}{2}} \|\tilde{G}\|_{L^2_{\tau,\zeta}} \|\hat{f}\|_{L^2_\zeta} \|\hat{g}\|_{L^2_\zeta}.
		\end{multline}
		By renaming the components, we can assume that $|\xi| \sim |\xi_1|\sim M_1$ and $|\eta| \sim |\eta_1| \sim M_2$, where $\xi = (\xi_1, \underline{\xi}), \eta= (\eta_1, \underline{\eta})$ with $\underline{\xi}, \underline{\eta} \in \R^{N-1}$. By the change of variables 
		\[
		\left\{
		\begin{array}{rcl}
		\tau &=& |\xi|^4 + |\eta|^4 - \mu (|\xi|^2 + |\eta|^2), \\
		\zeta &=& \xi + \eta,
		\end{array}
		\right.
		\]
		a direct computation shows
		\[
		d\tau d\zeta = J d\xi_1 d\eta,
		\]
		where
		\[
		J= |4(|\xi|^2 \xi_1 - |\eta|^2 \eta_1) - 2\mu(\xi_1-\eta_1)| \sim |\eta|^3 \sim M_2^3.
		\]
		By the Cauchy-Schwarz inequality with the fact $|\underline{\xi}| \lesssim M_1$, we get
		\begin{align*}
		\text{LHS} \eqref{bili-est-proof} &= \Big| \iiint_{\R \times \R^{N-1} \times \R^N} \overline{\tilde{G}}(\tau, \zeta) \widehat{P_{M_1}f}(\xi) \widehat{P_{M_2}g}(\eta) J^{-1} d\tau d\underline{\xi} d\zeta \Big| \\
		&\leq \|\tilde{G}\|_{L^2_{\tau,\zeta}} \int_{\R^{N-1}} \Big( \iint_{\R \times \R^N} |\widehat{P_{M_1} f}(\xi)|^2 |\widehat{P_{M_2}g}(\eta)|^2 J^{-2} d\tau d\zeta \Big)^{1/2} d\underline{\xi} \\
		&\lesssim M_1^{\frac{N-1}{2}}\|\tilde{G}\|_{L^2_{\tau,\zeta}}  \Big( \iiint_{\R\times \R^{N-1} \times\R^N} |\widehat{P_{M_1}f}(\xi)|^2 |\widehat{P_{M_2}g}(\eta)|^2 J^{-2} d\tau d\underline{\xi} d\zeta \Big)^{1/2} \\
		&\lesssim M_1^{\frac{N-1}{2}} \|\tilde{G}\|_{L^2_{\tau,\zeta}} \Big(\iint_{\R\times \R^{N-1} \times \R^N} |\widehat{P_{M_1} f}(\xi)|^2 |\widehat{P_{M_2} g}(\eta)|^2 J^{-1} d\xi_1 d\underline{\xi} d\eta \Big)^{1/2} \\
		&\leq M_1^{\frac{N-1}{2}} M_2^{-\frac{3}{2}} \|\tilde{G}\|_{L^2_{\tau,\zeta}} \|\widehat{P_{M_1} f}\|_{L^2_\zeta} \|\widehat{P_{M_2} g}\|_{L^2_\zeta}
		\end{align*}
		which proves \eqref{bili-est-proof}, and the proof is complete.
	\end{proof}
	
	We have the following bilinear estimate adapted to the $X^\gamma$- and $Y^\gamma$-spaces.
	\begin{lemma} \label{lem-bili-est-Y}
		Let $N\geq 5$ and $\mu \geq 0$. Then for any $0<T \leq \infty$ and $M_1 \leq M_2$, we have
		\begin{align}  \label{bili-est-Y}
		\|P_{M_1}u_1 P_{M_2} u_2\|_{L^2_{t,x}([0,T)\times \R^N)} &\lesssim M_1^{\frac{N-4}{2}} \left(\frac{M_1}{M_2}\right)^{\frac{3}{2}-} \|P_{M_1} u_1\|_{Y^0([0,T))} \|P_{M_2} u_2\|_{Y^0([0,T))}. 
		\end{align}
	\end{lemma}
	
	\begin{remark}
		By \eqref{embed}, we can replace the $Y^0([0,T))$-norm of the above estimate by the $X^0([0,T))$-norm.
	\end{remark}

	\noindent {\it Proof of Lemma $\ref{lem-bili-est-Y}$.}
		We follow the argument of \cite{BOP-TAMS}. 
		By \eqref{bili-est} and the tranference principle, we have
		\begin{align} \label{bili-est-Y-proof-1}
		\|P_{M_1} u_1 P_{M_2} u_2\|_{L^2_{t,x}([0,T)\times \R^N)} \lesssim M_1^{\frac{N-4}{2}} \left(\frac{M_1}{M_2}\right)^{\frac{3}{2}} \|P_{M_1} u_1\|_{U^2_\Delta([0,T),L^2(\R^N))} \|P_{M_2} u_2\|_{U^2_\Delta([0,T),L^2(\R^N))}.
		\end{align}
		Since $\left(4,\frac{2N}{N-2}\right)$ is Biharmonic admissible, we have
		\[
		\|e^{it\Delta} P_{M_1} f\|_{L^4_{t,x}([0,T)\times \R^N)} \lesssim \||\nabla|^{\frac{N-4}{4}} e^{it\Delta} P_{M_1}f\|_{L^4_tL^{\frac{2N}{N-2}}_x([0,T)\times \R^N)} \lesssim M_1^{\frac{N-4}{4}} \|f\|_{L^2(\R^N)}.
		\]
		The tranference principle gives
		\[
		\|P_{M_1} u\|_{L^4_{t,x}([0,T)\times \R^N)} \lesssim M_1^{\frac{N-4}{4}} \|u\|_{U^4_\Delta([0,T),L^2(\R^N))}.
		\]
		By H\"older's inequality, we get
		\begin{align} \label{bili-est-Y-proof-2}
		\|P_{M_1} u_1 P_{M_2} u_2\|_{L^2_{t,x}([0,T)\times \R^N)} \lesssim M_1^{\frac{N-4}{4}} M_2^{\frac{N-4}{4}} \|u_1\|_{U^4_\Delta ([0,T),L^2(\R^N))} \|u_2\|_{U^4_\Delta([0,T),L^2(\R^N))}.
		\end{align}
		By the interpolation lemma, we have from \eqref{bili-est-Y-proof-1} and \eqref{bili-est-Y-proof-2} that
		\begin{align*}
		\|P_{M_1} u_1 P_{M_2} u_2\|&_{L^2_{t,x}([0,T)\times \R^N)} \\
		&\lesssim M_1^{\frac{N-4}{2}} \left(\frac{M_1}{M_2}\right)^{\frac{3}{2}} \left( \ln \left[ \left(\frac{M_2}{M_1}\right)^{\frac{N+2}{4}}\right] +1\right)^2 \|u_1\|_{V^2_\Delta([0,T),L^2(\R^N))} \|u_2\|_{V^2_\Delta([0,T),L^2(\R^N))} \\
		&\lesssim   M_1^{\frac{N-4}{2}} \left(\frac{M_1}{M_2}\right)^{\frac{3}{2}-} \|u_1\|_{Y^0([0,T))} \|u_2\|_{Y^0([0,T))}.
		\end{align*}
		This shows \eqref{bili-est-Y} and the proof is complete.
	\hfill $\Box$
	
	To finish this section, we recall the following linear estimates which are needed in the sequel.
	\begin{lemma}[Linear estimates {\cite[Propositions 2.10 and 2.11]{HTT}}] \label{lem-lin-est}
		Let $N\geq 1$ and $\mu \geq 0$. Let $\gamma \geq 0$ and $0 <T\leq \infty$. Then it holds that 
		\begin{align}  \label{homo-est}
		\|U_\mu(t) f\|_{X^\gamma([0,T))} \leq \|f\|_{H^\gamma(\R^N)}
		\end{align}
		and
		\begin{align} \label{inhomo-est}
		\left\| \int_0^t U_\mu(t-s) F(s) ds\right\|_{X^\gamma([0,T))} \leq \sup_{v \in Y^{-\gamma}([0,T)) \atop \|v\|_{Y^{-\gamma}} =1} \left| \int_0^T \int_{\R^N} F(t,x) \overline{v}(t,x) dxdt \right|
		\end{align}
		for all $f\in H^\gamma(\R^N)$ and $F \in L^1([0,T), H^\gamma(\R^N))$.
	\end{lemma}

\section{Probabilistic nonlinear estimates}
\label{S3}
\setcounter{equation}{0}
In this section, we will prove probabilistic nonlinear estimates needed to show the almost sure well-posedness. Denote
\begin{align} \label{defi-Phi-v}
\Phi(v(t)) := \mp i \int_0^t U_\mu(t-s) \Nc(v+z)(s) ds
\end{align}
and
\begin{align} \label{defi-tilde-Phi-v}
\tilde{\Phi}(v(t)):= \mp i \int_0^t U_\mu(t-s) \Nc(v+\vareps z) (s) ds,
\end{align}
where $z, v$ are as in \eqref{z-omega} and $\Nc(f) = f \overline{f} f$. We have the following probabilistic nonlinear estimates.
\begin{proposition} \label{prop-proba-non-est}
	Let $N\geq 5$, $\mu \geq 0$ and $\gamma \in (\gamma_N, \gamc)$, where $\gamc$ and $\gamma_N$ are as in \eqref{defi-gamc} and \eqref{defi-gam-N} respectively. Let $f \in H^\gamma(\R^N)$ and $f^\omega$ be the Wiener randomization defined in \eqref{defi-random} satisfying \eqref{cond-distri}. 
	\begin{itemize}
		\item Let $0 < T \leq 1$. Then there exists $0 <\vartheta \ll 1$ such that
		\begin{align} 
		\|\Phi (v)\|_{X^{\frac{N-4}{2}}([0,T))} &\leq C_1 \left( \|v\|^3_{X^{\frac{N-4}{2}}([0,T))} + T^{\vartheta} R^3\right), \label{proba-non-est-1} \\
		\|\Phi(v_1)-\Phi(v_2)\|_{X^{\frac{N-4}{2}}([0,T))} &\leq C_2\left(\sum_{j=1}^2 \|v_j\|^2_{X^{\frac{N-4}{2}}([0,T))} + T^\vartheta R^2 \right) \|v_1-v_2\|_{X^{\frac{N-4}{2}}([0,T))} \label{proba-non-est-2}
		\end{align}
		for all $v,v_1,v_2 \in X^{\frac{N-4}{2}}([0,T))$ and all $R>0$, outside a set of probability at most
		\[
		C \exp \left(-c R^2 \|f\|^{-2}_{H^\gamma(\R^N)}\right).
		\]
		\item For $0<\vareps \ll 1$, we have
		\begin{align}
		\|\tilde{\Phi}(v)\|_{X^{\frac{N-4}{2}}(\R)} &\leq C_3 \left( \|v\|^3_{X^{\frac{N-4}{2}}(\R)} + R^3\right), \label{global-proba-non-est-1} \\
		\|\tilde{\Phi}(v_1)-\tilde{\Phi}(v_2)\|_{X^{\frac{N-4}{2}}(\R)} &\leq C_4 \left( \sum_{j=1}^2 \|v_j\|^2_{X^{\frac{N-4}{2}}(\R)} + R^2 \right) \|v_1-v_2\|_{X^{\frac{N-4}{2}}(\R)} \label{global-proba-non-est-2}
		\end{align}
		for all $v,v_1,v_2 \in X^{\frac{N-4}{2}}(\R)$ and all $R>0$, outside a set of probability at most
		\[
		C \exp \left(-cR^2 \vareps^{-2} \|f\|^{-2}_{H^\gamma(\R^N)}\right).
		\]
	\end{itemize}
\end{proposition}

\begin{proof}
	We mainly follow the argument of B\'enyi-Oh-Pocovnicu \cite{BOP-TAMS}. 
	
$\bullet$ Let $0<T \leq 1$. We only prove \eqref{proba-non-est-1}, and the one for \eqref{proba-non-est-2} is treated similarly.  Given $M\geq 1$, we define
\[
\Phi_M(v(t)) := \mp i \int_0^t U_\mu(t-s) P_{\leq M} \Nc(v+z)(s) ds.
\]
By Bernstein's and H\"older inequalities, we have
\begin{align}
\|P_{\leq M} \Nc(v+z)\|_{L^1_t H^{\frac{N-4}{2}}_x([0,T)\times \R^N)} &\lesssim M^{\frac{N-4}{2}} \|\Nc(v+z)\|_{L^1_t L^2_x([0,T)\times \R^N)} \nonumber \\
&\lesssim M^{\frac{N-4}{2}} \|v\|^3_{L^3_t L^6_x([0,T) \times \R^N)} + M^{\frac{N-4}{2}} \|z\|^3_{L^3_t L^6_x([0,T) \times\R^N)}. \label{est-PM}
\end{align}
Thanks to the local in time probabilistic Strichartz estimates, the second term in \eqref{est-PM} is finite almost surely. On the other hand, by the Sobolev embedding and \eqref{est-LqLr-Y}, we have
\[
\|v\|_{L^3_t L^6_x([0,T)\times \R^N)} \lesssim \||\nabla|^{\frac{N-4}{3}} v\|_{L^3_t L^{\frac{6N}{3N-8}}_x([0,T)\times \R^N)} \lesssim \|v\|_{X^{\frac{N-4}{2}}([0,T))} <\infty.
\]
This shows that for each $M\geq 1$, $P_{\leq M} \Nc(v+z) \in L^1_t H^{\frac{N-4}{2}}_x([0,T)\times \R^N)$ almost surely. Thus, by Lemma $\ref{lem-lin-est}$, 
\begin{align} \label{est-Phi-M}
\|\Phi_M(v)\|_{X^{\frac{N-4}{2}}([0,T))} \lesssim \sup_{v_4 \in Y^0([0,T)) \atop \|v_4\|_{Y^0}=1}  \left| \int_0^T \int_{\R^N} \scal{\nabla}^{\frac{N-4}{2}} [\Nc(v+z)(t,x)] \overline{v}_4(t,x) dx dt\right|,
\end{align}
where $v_4= P_{\leq M} v_4$. In the following, we estimate the right hand side of \eqref{est-Phi-M} independently of the cutoff size $M\geq 1$, by performing a case-by-case analysis of expressions of the form
\begin{align} \label{proba-non-est-proof-1}
\left| \int_0^T \int_{\R^N} \scal{\nabla}^{\frac{N-4}{2}} (w_1 w_2 w_3) v_4 dx dt \right|,
\end{align}
where $\|v_4\|_{Y^0([0,T))} \leq 1$ and $w_j = v$ or $z$, $j=1,2,3$. As a result, by letting $M \rightarrow \infty$, the same estimate holds for $\Phi(v)$ without any cutoff, thus yielding \eqref{proba-non-est-1}.

Before proceeding further, let us simplify some of the notation. In the following, we drop the complex conjugate sign since it plays no role. We also denote $X^\gamma([0,T))$ and $Y^0([0,T))$ by $X^\gamma$ and $Y^0$ since $T$ is fixed. Similarly, we will use $L^p_t L^q_x$, $L^p_{t,x}$ and $H^\gamma$ instead of $L^p_tL^q_x([0,T)\times \R^N)$, $L^p_{t,x}([0,T) \times \R^N)$ and $H^\gamma(\R^N)$ respectively. Lastly, in most of the cases, we dyadically decompose $w_j= v_j$ or $z_j$, $j=1,2,3$, and $v_4$ such that their spatial frequency supports are $\{|\xi_j| \sim M_j\}$ for some dyadic $M_j \geq 1$ but still denoted them as $w_j=v_j$ or $z_j$, $j=1,2,3$, and $v_4$. Note that by the Parseval formula \footnote{We have
	\[
	\int_{\R^N} f_1(x) f_2(x) f_3(x) f_4(x) dx = \int_{\xi_1+\xi_2+\xi_3+\xi_4=0} \hat{f}_1(\xi_1) \hat{f}_2(\xi_2) \hat{f}_3(\xi_3) \hat{f}_4(\xi_4),
	\]
	where $\mathlarger{\int}_{\xi_1+\xi_2+\xi_3+\xi_4}$ denotes the integration with respect to the hyperplane's measure $\delta_0(\xi_1+\xi_2+\xi_3+\xi_4) d\xi_1 d\xi_2 d\xi_3 d\xi_4$.}, we have $\sum_{j=1}^4 \xi_j =0$.

{\bf Case 1. $vvv$ case.} In this case, we do not need to perform dyadic decomposition, and we divide the frequency spaces into $\{|\xi_1| \geq |\xi_2|, |\xi_3|\}$, $\{|\xi_2| \geq |\xi_1|, |\xi_3|\}$, and $\{|\xi_3| \geq |\xi_1|, |\xi_2|\}$. Without loss of generality, we assume that $|\xi_1| \geq |\xi_2|, |\xi_3|$. By H\"older's inequality, we have
\begin{align*}
\left| \int_0^T \int_{\R^N} \scal{\nabla}^{\frac{N-4}{2}} v_1 v_2 v_3 v_4 dx dt\right| \leq \|\scal{\nabla}^{\frac{N-4}{2}} v_1\|_{L^{\frac{2(N+4)}{N}}_{t,x}} \|v_2\|_{L^{\frac{N+4}{2}}_{t,x}} \|v_3\|_{L^{\frac{N+4}{2}}_{t,x}} \|v_4\|_{L^{\frac{2(N+4)}{N}}_{t,x}}.
\end{align*}
By \eqref{est-Lp-Y} and Remark $\ref{rem-est-Y}$, we have
\[
\|\scal{\nabla}^{\frac{N-4}{2}} v_1\|_{L^{\frac{2(N+4)}{N}}_{t,x}} \lesssim \|v_1\|_{X^{\frac{N-4}{2}}}, \quad \|v_4\|_{L^{\frac{2(N+4)}{N}}_{t,x}} \lesssim \|v_4\|_{Y^0} \leq 1.
\]
By Sobolev embedding and \eqref{est-LqLr-Y}, we have
\[
\|v_2\|_{L^{\frac{N+4}{2}}_{t,x}} \lesssim \||\nabla|^{\frac{N-4}{2}} v_2\|_{L^{\frac{N+4}{2}}_t L^{\frac{2N(N+4)}{N^2+4N-16}}_x} \lesssim \|v_2\|_{X^{\frac{N-4}{2}}},
\]
similarly for $v_3$. We thus get
\[
\left| \int_0^T \int_{\R^N} \scal{\nabla}^{\frac{N-4}{2}} v_1 v_2 v_3 v_4 dx dt\right| \lesssim \prod_{j=1}^3 \|v_j\|_{X^{\frac{N-4}{2}}}.
\]

{\bf Case 2. $zzz$ case.} Without loss of generality, we assume $M_3 \geq M_2 \geq M_1$. Note that $M_4 \lesssim M_3$ since otherwise the corresponding localized functions have disjoint supports.

{\bf Subcase 2a. $M_2 \sim M_3$.} By H\"older's inequality, we have
\begin{align*}
\left| \int_0^T \int_{\R^N} z_1 z_2\scal{\nabla}^{\frac{N-4}{2}} z_3 v_4 dx dt\right| &\leq \|z_1\|_{L^{\frac{N+4}{2}}_{t,x}} \|z_2\|_{L^4_{t,x}} \|\scal{\nabla}^{\frac{N-4}{2}} z_3\|_{L^4_{t,x}} \|v_4\|_{L^{\frac{2(N+4)}{N}}_{t,x}} \\
&\sim \|z_1\|_{L^{\frac{N+4}{2}}_{t,x}} \|\scal{\nabla}^{\frac{N-4}{4}} z_2\|_{L^4_{t,x}} \|\scal{\nabla}^{\frac{N-4}{4}} z_3\|_{L^4_{t,x}} \|v_4\|_{L^{\frac{2(N+4)}{N}}_{t,x}}.
\end{align*}
By the Littlewood-Paley decomposition and the local-in-time probabilistic Strichartz estimates, we observe that for $r \geq 2$,
\begin{align*}
\sum_{M_1 \geq 1} \|P_{M_1} z_1\|_{L^r_{t,x}} &= \sum_{M_1 \geq 1} M_1^{0-} \|\scal{\nabla}^{0+} P_{M_1} z_1\|_{L^r_{t,x}} \\
&\leq \Big(\sum_{M_1\geq 1} M_1^{0-}\Big)^{\frac{1}{r'}} \Big\| \|\scal{\nabla}^{0+} P_{M_1} z_1\|_{L^r_{t,x}} \Big\|_{\ell^r_{M_1}} \\
&= \Big(\sum_{M_1\geq 1} M_1^{0-}\Big)^{\frac{1}{r'}} \Big\| \|\scal{\nabla}^{0+} P_{M_1} z_1\|_{\ell^r_{M_1}} \Big\|_{L^r_{t,x}} \\
&\leq \Big(\sum_{M_1\geq 1} M_1^{0-}\Big)^{\frac{1}{r'}} \Big\| \|\scal{\nabla}^{0+} P_{M_1} z_1\|_{\ell^2_{M_1}} \Big\|_{L^r_{t,x}} \\
&\lesssim \|\scal{\nabla}^{0+} z_1\|_{L^r_{t,x}} \\
&= \|\scal{\nabla}^{0+} U_\mu(t) f^\omega\|_{L^r_{t,x}} \\
&\leq T^{0+} R
\end{align*}
outside a set of probability at most 
\[
C \exp \left(-c R^2T^{-\frac{2}{r}+} \|f\|^{-2}_{H^{0+}}\right).
\]
Applying the above observation to $r=\frac{N+4}{2}$ and $4$, we obtain
\[
\left| \int_0^T \int_{\R^N} z_1 z_2\scal{\nabla}^{\frac{N-4}{2}} z_3 v_4 dx dt\right| \lesssim T^{0+}R^3
\]
outside a set of probability at most
\[
C \exp \left(-cR^2T^{-\frac{4}{N+4}+} \|f\|^{-2}_{H^{0+}}\right) + C \exp \left(-cR^2T^{-\frac{1}{2}+} \|f\|^{-2}_{H^{\frac{N-4}{4}+}}\right).
\]
Here we have extracted a negative power of $M_3$ and used $M_4\lesssim M_3$ to estimate
\[
\sum_{M_4\geq 1} M_4^{0-} \|P_{M_4} v_4\|_{L^{\frac{2(N+4)}{N}}_{t,x}} \lesssim \|v_4\|_{L^{\frac{2(N+4)}{N}}_{t,x}} \lesssim \|v_4\|_{Y^0} \leq 1.
\]
Note that since $T\leq 1$ and $\gamma >\frac{N-4}{4}$, the above probability can be bounded by
\[
C \exp \left(-cR^2\|f\|^{-2}_{H^\gamma} \right).
\]

{\bf Subcase 2b. $M_3 \gg M_2 \geq M_1$.} Note that we must have $M_4 \sim M_3$. 

{\bf Subcase 2b.i. $M_3^{\frac{3}{N-1}} \gg M_2 \geq M_1$.} By H\"older's inequality, we have
\begin{align*}
\left| \int_0^T \int_{\R^N} z_1 z_2 \scal{\nabla}^{\frac{N-4}{2}} z_3 v_4 dx dt \right| \leq \|z_2 \scal{\nabla}^{\frac{N-4}{2}} z_3\|_{L^2_{t,x}} \|z_1 v_4\|_{L^2_{t,x}}.
\end{align*}
Let $a>0$ be a small constant to be chosen shortly, we estimate
\begin{align*}
\|z_2 \scal{\nabla}^{\frac{N-4}{2}} z_3\|_{L^2_{t,x}} &\lesssim M_3^{\frac{(N-4)a}{2}}\|z_2\|_{L^4_{t,x}}^a \|z_3\|_{L^4_{t,x}}^a \|z_2 \scal{\nabla}^{\frac{N-4}{2}} z_3\|^{1-a}_{L^2_{t,x}} \\
&\lesssim M_2^{-a\gamma} M_3^{\frac{(N-4)a}{2}-a\gamma} \|\scal{\nabla}^\gamma z_2\|^a_{L^4_{t,x}} \|\scal{\nabla}^\gamma z_3\|^a_{L^4_{t,x}} \|z_2 \scal{\nabla}^{\frac{N-4}{2}} z_3\|^{1-a}_{L^2_{t,x}}.
\end{align*}
We next use \eqref{bili-est-Y} and \eqref{prop-X-gamma} to have
\begin{align*}
\|z_2 \scal{\nabla}^{\frac{N-4}{2}} z_3\|_{L^2_{t,x}} &\lesssim M_2^{\frac{N-1}{2}-} M_3^{-\frac{3}{2}+} \|z_2\|_{X^0} \|\scal{\nabla}^{\frac{N-4}{2}} z_3\|_{X^0} \\
&\lesssim M_2^{\frac{N-1}{2}-} M_3^{\frac{N-4}{2}-\frac{3}{2}+} \|P_{M_2} f^\omega\|_{L^2} \|P_{M_3}f^\omega\|_{L^2} \\
&\lesssim M_2^{\frac{N-1}{2}-\gamma-} M_3^{\frac{N-4}{2}-\gamma-\frac{3}{2}+} \|P_{M_2}f^\omega\|_{H^\gamma} \|P_{M_3} f^\omega\|_{H^\gamma}.
\end{align*}
It follows that
\begin{align} \label{est-a}
\|z_2 \scal{\nabla}^{\frac{N-4}{2}}z_3\|_{L^2_{t,x}} &\lesssim M_2^{\frac{N-1}{2} -\gamma - \frac{N-1}{2}a-} M_3^{\frac{N-7}{2} -\gamma +\frac{3a}{2}+} \prod_{j=2}^3 \|\scal{\nabla}^\gamma z_j\|^a_{L^4_{t,x}} \|P_{M_j}f^\omega\|_{H^\gamma}^{1-a}.
\end{align}
Similarly, by \eqref{bili-est-Y}, we have
\begin{align*}
\|z_1 v_4\|_{L^2_{t,x}} &\lesssim M_1^{\frac{N-1}{2}-} M_4^{-\frac{3}{2}+} \|z_1\|_{X^0} \|v_4\|_{Y^0} \\
&\lesssim M_1^{\frac{N-1}{2}-} M_4^{-\frac{3}{2}+} \|P_{M_1} f^\omega\|_{L^2} \\
&\lesssim M_1^{\frac{N-1}{2}-\gamma-} M_4^{-\frac{3}{2}+} \|P_{M_1}f^\omega\|_{H^\gamma}.
\end{align*}
Since $M_3^{\frac{3}{N-1}} \gg M_2 \geq M_1$ and $M_3 \sim M_4$, we obtain
\begin{align*}
\left| \int_0^T \int_{\R^N} z_1 z_2 \scal{\nabla}^{\frac{N-4}{2}} z_3 v_4 dx dt \right| &\lesssim M_3^{\frac{N-4}{2}-\frac{N+5}{N-1}\gamma+} \|P_{M_1}f^\omega\|_{H^\gamma} \prod_{j=2}^3 \|\scal{\nabla}^\gamma z_j\|^a_{L^4_{t,x}} \|P_{M_j}f^\omega\|_{H^\gamma}^{1-a}.
\end{align*}
provided
\[
\frac{N-1}{2}-\gamma -\frac{N-1}{2} a>0 \quad \text{or} \quad a<1-\frac{2\gamma}{N-1}.
\]
We want the largest frequency $M_3$ to have a negative power so that we can sum over dyadic blocks. This requires
\begin{align} \label{cond-gamma-proof}
\frac{N-4}{2}-\frac{N+5}{N-1}\gamma <0
\end{align}
which is satisfied as
\[
\gamma > \gamma_N \geq \frac{(N-1)(N-4)}{2(N+5)}.
\]
Under this condition, we can sum over dyadic blocks as in Case 2a. We thus get
\[
\left| \int_0^T \int_{\R^N} z_1 z_2 \scal{\nabla}^{\frac{N-4}{2}} z_3 v_4 dx dt \right| \lesssim T^{0+}R^3
\]
outside a set of probability at most
\[
C \exp \left(-cR^2T^{-\frac{1}{2}+} \|f\|^{-2}_{H^\gamma}\right) + C\exp\left(-cR^2\|f\|^{-2}_{H^\gamma}\right).
\]

{\bf Subcase 2b.ii. $M_2 \gtrsim M_3^{\frac{3}{N-1}} \gg M_1$.} By H\"older's inequality and \eqref{bili-est-Y}, we have
\begin{align*}
\Big| \int_0^T \int_{\R^N} z_1 z_2 &\scal{\nabla}^{\frac{N-4}{2}} z_3 v_4 dx dt \Big| \leq \|z_2\|_{L^4_{t,x}} \|\scal{\nabla}^{\frac{N-4}{2}} z_3\|_{L^4_{t,x}} \|z_1 v_4\|_{L^2_{t,x}} \\
&\lesssim M_1^{\frac{N-1}{2} - \gamma-} M_2^{-\gamma} M_3^{\frac{N-4}{2}-\gamma} M_4^{-\frac{3}{2}+} \|P_{M_1} f^\omega\|_{H^\gamma} \|\scal{\nabla}^\gamma z_2\|_{L^4_{t,x}} \|\scal{\nabla}^\gamma z_3\|_{L^4_{t,x}} \|v_4\|_{Y^0} \\
&\lesssim M_3^{\frac{N-4}{2}-\frac{N+5}{N-1}\gamma+} \|P_{M_1}f^\omega\|_{H^\gamma} \|\scal{\nabla}^\gamma z_2\|_{L^4_{t,x}} \|\scal{\nabla}^\gamma z_3\|_{L^4_{t,x}} \\
&\lesssim T^{0+}R^3
\end{align*}
outside a set of probability at most 
\[
C \exp \left(-c R^2T^{-\frac{1}{2}+} \|f\|^{-2}_{H^\gamma}\right) + C \exp \left(-cR^2\|f\|^{-2}_{H^\gamma}\right)
\]
as long as \eqref{cond-gamma-proof} holds.

{\bf Subcase 2b.iii. $M_2 \geq M_1 \gtrsim M_3^{\frac{3}{N-1}}$.} By H\"older's inequality, we have
\begin{align*}
\Big| \int_0^T \int_{\R^N} z_1 z_3 \scal{\nabla}^{\frac{N-4}{2}} z_3 v_4 dx dt \Big| &\lesssim \|z_1\|_{L^{\frac{6(N+4)}{N+8}}_{t,x}} \|z_2\|_{L^{\frac{6(N+4)}{N+8}}_{t,x}} \|\scal{\nabla}^{\frac{N-4}{2}} z_3\|_{L^{\frac{6(N+4)}{N+8}}_{t,x}} \|v_4\|_{L^{\frac{2(N+4)}{N}}_{t,x}} \\
&\lesssim M_1^{-\gamma} M_2^{-\gamma} M_3^{\frac{N-4}{2}-\gamma} \prod_{j=1}^3 \|\scal{\nabla}^\gamma z_j\|_{L^{\frac{6(N+4)}{N+8}}_{t,x}} \|v_4\|_{Y^0} \\
&\lesssim M_3^{\frac{N-4}{2}-\frac{N+5}{N-1}\gamma} \prod_{j=1}^3 \|\scal{\nabla}^\gamma z_j\|_{L^{\frac{6(N+4)}{N+8}}_{t,x}} \\
&\lesssim T^{0+} R^3
\end{align*}
outside a set of probability at most
\[
C \exp \left(-cR^2T^{-\frac{N+8}{3(N+4)}+} \|f\|^{-2}_{H^\gamma}\right)
\]
as long as \eqref{cond-gamma-proof} holds. 

{\bf Case 3. $vvz$ case.} Without loss of generality, we assume $M_1 \geq M_2$.

{\bf Subcase 3a. $M_1 \gtrsim M_3$.} In this case, we only perform the dyadic decomposition on $v_1, v_2$ and $z_3$. Note that $M_1 \gtrsim \max\{M_2, M_3, |\xi_4|\}$, where $\xi_4$ is the spatial frequency of $v_4$. By H\"older's inequality and \eqref{bili-est-Y}, we have
\begin{align*}
\Big| \int_0^T \int_{\R^N} \scal{\nabla}^{\frac{N-4}{2}} &v_1 v_2 z_3 v_4 dx dt \Big| \lesssim \sum_{M_1 \geq M_2} \|\scal{\nabla}^{\frac{N-4}{2}} P_{M_1}v_1 P_{M_2} v_2\|_{L^2_{t,x}} \sum_{M_3} \|P_{M_3} z_3\|_{L^{\frac{N+4}{2}}_{t,x}} \|v_4\|_{L^{\frac{2(N+4)}{N}}_{t,x}} \\
&\lesssim \sum_{M_1\geq M_2} M_2^{\frac{N-4}{2}} \left(\frac{M_2}{M_1}\right)^{\frac{3}{2}-} \|\scal{\nabla}^{\frac{N-4}{2}} P_{M_1}v_1\|_{X^0} \|P_{M_2}v_2\|_{X^0} \sum_{M_3}\|P_{M_3}z_3\|_{L^{\frac{N+4}{2}}_{t,x}} \|v_4\|_{Y^0} \\
&\lesssim \sum_{M_1\geq M_2} \left(\frac{M_2}{M_1}\right)^{\frac{3}{2}-} \|P_{M_1}v_1\|_{X^{\frac{N-4}{2}}} \|P_{M_2}v_2\|_{X^{\frac{N-4}{2}}} \sum_{M_3}\|P_{M_3}z_3\|_{L^{\frac{N+4}{2}}_{t,x}}.
\end{align*}
If $M_1 \sim M_2$, then we simply remove $\left(\frac{M_2}{M_1}\right)^{\frac{3}{2}-}$ and use Cauchy-Schwarz inequality to bound
	\begin{align*}
	\sum_{M_1 \sim M_2} \|P_{M_1}v_1\|_{X^{\frac{N-4}{2}}} \|P_{M_2}v_2\|_{X^{\frac{N-4}{2}}} &\lesssim \left(\sum_{M1} \|P_{M_1}v_1\|^2_{X^{\frac{N-4}{2}}}\right)^{\frac{1}{2}} \left(\sum_{M_1} \|P_{M_1}v_2\|^2_{X^{\frac{N-4}{2}}}\right)^{\frac{1}{2}} \\
	&\lesssim \|v_1\|_{X^{\frac{N-4}{2}}} \|v_2\|_{X^{\frac{N-4}{2}}}.
	\end{align*}
	If $M_1 \gg M_2$, then we can extract from $\left(\frac{M_2}{M_1}\right)^{\frac{3}{2}-}$ a negative power of $M_1$ which allows to sum over $M_1$. By extracting a negative power of $M_1$, we can sum over $M_2$. We thus get
	\[
	\Big| \int_0^T \int_{\R^N} \scal{\nabla}^{\frac{N-4}{2}} v_1 v_2 z_3 v_4 dx dt \Big| \lesssim T^{0+} R \|v_1\|_{X^{\frac{N-4}{2}}} \|v_2\|_{X^{\frac{N-4}{2}}}
	\]
	outside a set of probability at most
	\[
	C \exp \left(-cR^2T^{-\frac{4}{N+4}+}\|f\|^{-2}_{H^{0+}}\right).
	\]

{\bf Subcase 3b. $M_3 \gg M_1 \geq M_2$.} Note that we must have $M_3 \sim M_4$. 

{\bf Subcase 3b.i. $M_1 \gtrsim M_3^{\frac{3}{N-1}}$.} By H\"older's inequality, Lemma $\ref{lem-str-est-Y}$ and \eqref{bili-est-Y}, we have
\begin{align*}
\Big|\int_0^T \int_{\R^N} v_1 v_2 \scal{\nabla}^{\frac{N-4}{2}} z_3 v_4 dx dt \Big| &\leq \|v_1\|_{L^{\frac{2(N+4)}{N}}_{t,x}} \|\scal{\nabla}^{\frac{N-4}{2}} z_3\|_{L^{\frac{N+4}{2}}_{t,x}} \|v_2 v_4\|_{L^2_{t,x}} \\
&\lesssim M_2^{\frac{N-4}{2}} \left(\frac{M_2}{M_4}\right)^{\frac{3}{2}-} M_3^{\frac{N-4}{2}} \|v_1\|_{X^0} \|v_2\|_{X^0} \|z_3\|_{L^{\frac{N+4}{2}}_{t,x}}  \|v_4\|_{Y^0} \\
&\lesssim M_1^{-\frac{N-4}{2}} M_2^{\frac{3}{2}-} M_3^{\frac{N-7}{2}-\gamma+} \|v_1\|_{X^{\frac{N-4}{2}}} \|v_2\|_{X^{\frac{N-4}{2}}} \|\scal{\nabla}^\gamma z_3\|_{L^{\frac{N+4}{2}}_{t,x}} \\
&\lesssim M_3^{\frac{(N-4)(N-7)}{2(N-1)} - \gamma +} \|v_1\|_{X^{\frac{N-4}{2}}} \|v_2\|_{X^{\frac{N-4}{2}}} \|\scal{\nabla}^\gamma z_3\|_{L^{\frac{N+4}{2}}_{t,x}} \\
&\lesssim T^{0+}R \|v_1\|_{X^{\frac{N-4}{2}}} \|v_2\|_{X^{\frac{N-4}{2}}}
\end{align*}
outside a set of probability at most
\[
C \exp \left(-cR^2T^{-\frac{4}{N+4}+} \|f\|^{-2}_{H^\gamma} \right)
\]
provided
\[
\frac{(N-4)(N-7)}{2(N-1)} - \gamma  <0 \quad \text{or} \quad \gamma>\frac{(N-4)(N-7)}{2(N-1)}
\]
which is less restrictive than \eqref{cond-gamma-proof}.

{\bf Subcase 3b.ii. $M_3^{\frac{3}{N-1}} \gg M_1$.} By H\"older's inequality, we have
\begin{align*}
\Big| \int_0^T \int_{\R^N} v_1 v_2 \scal{\nabla}^{\frac{N-4}{2}} z_3 v_4 dx dt \Big| \leq \|v_1 \scal{\nabla}^{\frac{N-4}{2}} z_3\|_{L^2_{t,x}} \|v_2v_4\|_{L^2_{t,x}}.
\end{align*}
Let $a>0$ be a small constant to be chosen later, we use Lemma $\ref{lem-str-est-Y}$ and \eqref{bili-est-Y} to estimate
\begin{align*}
\|v_1 \scal{\nabla}^{\frac{N-4}{2}} z_3\|_{L^2_{t,x}} &\lesssim M_3^{\frac{(N-4)a}{2}} \|v_1\|^a_{L^{\frac{2(N+4)}{N}}_{t,x}} \|z_3\|^a_{L^{\frac{N+4}{2}}_{t,x}} \|v_1 \scal{\nabla}^{\frac{N-4}{2}} z_3\|^{1-a}_{L^2_{t,x}} \\
&\lesssim M_1^{\frac{(N-1)(1-a)}{2}-} M_3^{-\frac{3(1-a)}{2}+} M_3^{\frac{(N-4)a}{2}} \|v_1\|_{X^0} \|z_3\|^a_{L^{\frac{N+4}{2}}_{t,x}} \|P_{M_3}\scal{\nabla}^{\frac{N-4}{2}}f^\omega\|^{1-a}_{L^2} \\
&\lesssim M_1^{\frac{3}{2} - \frac{(N-1)a}{2} -} M_3^{\frac{N-7}{2} - \gamma +\frac{3a}{2}+} \|v_1\|_{X^{\frac{N-4}{2}}} \|\scal{\nabla}^\gamma z_3\|^a_{L^{\frac{N+4}{2}}_{t,x}} \|P_{M_3}f^\omega\|^{1-a}_{H^\gamma}.
\end{align*}
Similarly, we have
\begin{align*}
\|v_2v_4\|_{L^2_{t,x}} &\lesssim M_2^{\frac{N-1}{2} -} M_4^{-\frac{3}{2}+} \|v_2\|_{X^0} \|v_4\|_{Y^0} \\
& \lesssim M_2^{\frac{3}{2}-} M_3^{-\frac{3}{2}+} \|v_2\|_{X^{\frac{N-4}{2}}}.
\end{align*}
It follows that
\begin{align*}
\Big| \int_0^T \int_{\R^N} v_1 v_2 &\scal{\nabla}^{\frac{N-4}{2}} z_3 v_4 dx dt \Big|\\ &\lesssim M_1^{\frac{3}{2}-\frac{(N-1)a}{2} -} M_2^{\frac{3}{2}-} M_3^{\frac{N-10}{2} -\gamma +\frac{3a}{2}+} \|v_1\|_{X^{\frac{N-4}{2}}} \|v_2\|_{X^{\frac{N-4}{2}}} \|\scal{\nabla}^\gamma z_3\|^a_{L^{\frac{N+4}{2}}_{t,x}} \|P_{M_3}f^\omega\|^{1-a}_{H^\gamma} \\
&\lesssim M_3^{\frac{(N-4)(N-7)}{2(N-1)} -\gamma+}  \|v_1\|_{X^{\frac{N-4}{2}}} \|v_2\|_{X^{\frac{N-4}{2}}} \|\scal{\nabla}^\gamma z_3\|^a_{L^{\frac{N+4}{2}}_{t,x}} \|P_{M_3}f^\omega\|^{1-a}_{H^\gamma}
\end{align*}
provided
\[
\frac{3}{2}-\frac{(N-1)a}{2} >0 \quad \text{or} \quad a<\frac{3}{N-1}.
\]
We want the power of $M_3$ is strictly negative in order to sum over dyadic blocks. This requires
\[
\gamma >\frac{(N-4)(N-7)}{2(N-1)}
\]
which is less restrictive than \eqref{cond-gamma-proof}. It follows that
\begin{align*}
\Big| \int_0^T \int_{\R^N} v_1 v_2 \scal{\nabla}^{\frac{N-4}{2}} z_3 v_4 dx dt \Big|
\lesssim T^{0+} R \|v_1\|_{X^{\frac{N-4}{2}}} \|v_2\|_{X^{\frac{N-4}{2}}}
\end{align*}
outside a set of probability at most
\[
C\exp \left(-cR^2T^{-\frac{4}{N+4}+} \|f\|^{-2}_{H^\gamma} \right) + C \exp \left(-cR^2\|f\|^{-2}_{H^\gamma} \right).
\]

{\bf Case 4. $vzz$ case.} Without loss of generality, we assume that $M_3 \geq M_2$. 

{\bf Subcase 4a. $M_1 \gtrsim M_3$.} By H\"older's inequality and Lemma $\ref{lem-str-est-Y}$, we have
\begin{align*}
\Big| \int_0^T \int_{\R^N} \scal{\nabla}^{\frac{N-4}{2}} v_1 z_2 z_3 v_4 dx dt \Big| &\leq \|\scal{\nabla}^{\frac{N-4}{2}} v_1\|_{L^{\frac{2(N+4)}{N}}_{t,x}} \|z_2\|_{L^{\frac{N+4}{2}}_{t,x}} \|z_3\|_{L^{\frac{N+4}{2}}_{t,x}} \|v_4\|_{L^{\frac{2(N+4)}{N}}_{t,x}} \\
&\lesssim \|v_1\|_{X^{\frac{N-4}{2}}} \|z_2\|_{L^{\frac{N+4}{2}}_{t,x}} \|z_3\|_{L^{\frac{N+4}{2}}_{t,x}} \|v_4\|_{Y^0} \\
&\lesssim T^{0+} R^2 \|v_1\|_{X^{\frac{N-4}{2}}}
\end{align*}
outside a set of probability at most
\[
C \exp \left(-cR^2T^{-\frac{4}{N+4}+} \|f\|^{-2}_{H^{0+}}\right).
\]
Here we have used that if $M_3 \gtrsim \max \{M_1, M_4\}$, then we can extract a negative power of $M_3$ to sum over $M_1$ and $M_4$. Otherwise, we have $M_1 \sim M_4  \gg M_3$. In this case, we can use Cauchy-Schwarz inequality to sum over $M_1$ and $M_4$. 

{\bf Subcase 4b. $M_3 \gg M_1$.} 

{\bf Subcase 4b.1. $M_3 \sim M_2 \gg M_1$.} We must have $M_1 \sim M_4$. By H\"older's inequality and Lemma $\ref{lem-str-est-Y}$, we have
\begin{align*}
\Big| \int_0^T \int_{\R^N} v_1 z_2\scal{\nabla}^{\frac{N-4}{2}} z_3 v_4 dx dt \Big| &\leq \|v_1\|_{L^{\frac{N+4}{2}}_{t,x}} \|z_2\|_{L^4_{t,x}} \|\scal{\nabla}^{\frac{N-4}{2}} z_3\|_{L^4_{t,x}} \|v_4\|_{L^{\frac{2(N+4)}{N}}_{t,x}} \\
&\lesssim M_3^{\frac{N-4}{2} - 2\gamma} \|v_1\|_{X^{\frac{N-4}{2}}} \|\scal{\nabla}^\gamma z_2\|_{L^4_{t,x}} \|\scal{\nabla}^\gamma z_3\|_{L^4_{t,x}} \|v_4\|_{Y^0} \\
&\lesssim T^{0+}R^2 \|v_1\|_{X^{\frac{N-4}{2}}} 
\end{align*}
outside a set of probability at most
\[
C \exp \left(-cR^2T^{-\frac{1}{2}+} \|f\|^{-2}_{H^\gamma} \right)
\]
as long as $\gamma>\frac{N-4}{4}$ which is less restrictive than \eqref{cond-gamma-proof}.

{\bf Subcase 4b.2. $M_3 \gg M_2, M_1$.} In this case, we must have $M_3 \sim M_4$.

{\bf Subcase 4b.2.i. $M_1, M_2 \ll M_3^{\frac{3}{N-1}}$.} By H\"older's inequality, \eqref{est-a} and Lemma $\ref{lem-str-est-Y}$, we have
\begin{align*}
\Big| \int_0^T \int_{\R^N} & v_1 z_2 \scal{\nabla}^{\frac{N-4}{2}} z_3 v_4 dx dt \Big| \leq \|z_2\scal{\nabla}^{\frac{N-4}{2}}z_3\|_{L^2_{t,x}} \|v_1 v_4\|_{L^2_{t,x}} \\
&\lesssim M_1^{\frac{N-1}{2}-} M_2^{\frac{N-1}{2}-\gamma-\frac{N-1}{2}a-} M_3^{\frac{N-7}{2}-\gamma+\frac{3a}{2}+} M_4^{-\frac{3}{2}+} \|v_1\|_{X^0} \prod_{j=2}^3 \left(\|\scal{\nabla}^\gamma z_j\|^a_{L^4_{t,x}} \|P_{M_j} f^\omega\|^{1-a}_{H^\gamma}\right) \|v_4\|_{Y^{0}} \\
&\lesssim M_1^{\frac{3}{2}-} M_2^{\frac{N-1}{2}-\gamma-\frac{N-1}{2}a-} M_3^{\frac{N-10}{2}-\gamma+\frac{3a}{2}+} \|v_1\|_{X^{\frac{N-4}{2}}} \prod_{j=2}^3 \left(\|\scal{\nabla}^\gamma z_j\|^a_{L^4_{t,x}} \|P_{M_j} f^\omega\|^{1-a}_{H^\gamma}\right) \\
&\lesssim M_3^{\frac{(N-4)^2}{2(N-1)} -\frac{(N+2)\gamma}{N-1} +} \|v_1\|_{X^{\frac{N-4}{2}}} \prod_{j=2}^3 \left(\|\scal{\nabla}^\gamma z_j\|^a_{L^4_{t,x}} \|P_{M_j} f^\omega\|^{1-a}_{H^\gamma}\right)
\end{align*}
provided 
\[
a<1-\frac{2\gamma}{N-1}.
\]
We want the largest frequency to have a negative power in order to sum over dyadic blocks. This requires 
\[
\gamma>\frac{(N-4)^2}{2(N+2)}
\]
which is again less restrictive than \eqref{cond-gamma-proof}. We thus get
\[
\Big| \int_0^T \int_{\R^N} v_1 z_2 \scal{\nabla}^{\frac{N-4}{2}} z_3 v_4 dx dt \Big| \lesssim T^{0+} R^2 \|v_1\|_{X^{\frac{N-4}{2}}}
\]
outside a set of probability at most
\[
C \exp \left(-cR^2T^{-\frac{1}{2}+} \|f\|^{-2}_{H^\gamma}\right) + C \exp \left(-cR^2\|f\|^{-2}_{H^\gamma}\right).
\]

{\bf Subcase 4b.2.ii. $M_1 \ll M_3^{\frac{3}{N-1}} \lesssim M_2$.} By H\"older's inequality, Lemma $\ref{lem-str-est-Y}$ and \eqref{bili-est-Y}, we have
\begin{align*}
\Big| \int_0^T \int_{\R^N} v_1 z_2 \scal{\nabla}^{\frac{N-4}{2}} z_3 v_4 dx dt \Big| &\leq \|z_2\|_{L^4_{t,x}} \|\scal{\nabla}^{\frac{N-4}{2}} z_3\|_{L^4_{t,x}} \|v_1 v_4\|_{L^2_{t,x}} \\
&\lesssim M_1^{\frac{N-1}{2}-} M_2^{-\gamma} M_3^{\frac{N-4}{2}-\gamma} M_4^{-\frac{3}{2}+} \|v_1\|_{X^0} \prod_{j=2}^3 \|\scal{\nabla}^\gamma z_j\|_{L^4_{t,x}} \|v_4\|_{Y^0} \\
&\lesssim M_1^{\frac{3}{2}-} M_2^{-\gamma} M_3^{\frac{N-7}{2}-\gamma+} \|v_1\|_{X^{\frac{N-4}{2}}} \prod_{j=2}^3 \|\scal{\nabla}^\gamma z_j\|_{L^4_{t,x}} \\
&\lesssim M_3^{\frac{(N-4)^2}{2(N-1)} -\frac{(N+2)\gamma}{N-1} +} \|v_1\|_{X^{\frac{N-4}{2}}} \prod_{j=2}^3 \|\scal{\nabla}^\gamma z_j\|_{L^4_{t,x}}\\
&\lesssim T^{0+} R^2 \|v_1\|_{X^{\frac{N-4}{2}}}
\end{align*}
outside a set of probability at most 
\[
C \exp \left(-cR^2T^{-\frac{1}{2}+} \|f\|^{-2}_{H^\gamma}\right)
\]
as long as 
\[
\gamma >\frac{(N-4)^2}{2(N+2)}
\]
which is satisfied if \eqref{cond-gamma-proof} holds.

{\bf Subcase 4b.2.iii. $M_2 \ll M_3^{\frac{3}{N-1}} \lesssim M_1$.} By H\"older's inequality, Lemma $\ref{lem-str-est-Y}$ and \eqref{bili-est-Y}, we have
\begin{align*}
\Big|\int_0^T \int_{\R^N} v_1 z_2 \scal{\nabla}^{\frac{N-4}{2}} z_3 v_4 dx dt \Big| &\leq \|v_1\|_{L^{\frac{2(N+4)}{N}}_{t,x}} \|\scal{\nabla}^{\frac{N-4}{2}} z_3\|_{L^{\frac{N+4}{2}}_{t,x}} \|z_2 z_4\|_{L^2_{t,x}} \\
&\lesssim M_2^{\frac{N-1}{2}-}M_3^{\frac{N-4}{2}-\gamma} M_4^{-\frac{3}{2}+} \|v_1\|_{X^0} \|P_{M_2}f^\omega\|_{L^2} \|\scal{\nabla}^\gamma z_3\|_{L^{\frac{N+4}{2}}_{t,x}} \|v_4\|_{Y^0} \\
&\lesssim M_1^{-\frac{N-4}{2}} M_2^{\frac{N-1}{2}-\gamma-} M_3^{\frac{N-7}{2}-\gamma+} \|v_1\|_{X^{\frac{N-4}{2}}} \|P_{M_2}f^\omega\|_{H^\gamma} \|\scal{\nabla}^\gamma z_3\|_{L^{\frac{N+4}{2}}_{t,x}} \\
&\lesssim M_3^{\frac{(N-4)^2}{2(N-1)} -\frac{(N+2)\gamma}{N-1} +}\|v_1\|_{X^{\frac{N-4}{2}}} \|P_{M_2}f^\omega\|_{H^\gamma} \|\scal{\nabla}^\gamma z_3\|_{L^{\frac{N+4}{2}}_{t,x}} \\
&\lesssim T^{0+} R^2 \|v_1\|_{X^{\frac{N-4}{2}}}
\end{align*}
outside a set of probability at most
\[
C \exp \left(-cR^2\|f\|^{-2}_{H^\gamma}\right) + C \exp \left(-cR^2T^{-\frac{4}{N+4}+} \|f\|^{-2}_{H^\gamma}\right)
\]
as long as 
\[
\gamma>\frac{(N-4)^2}{2(N+2)}
\]
which is satisfied if \eqref{cond-gamma-proof} holds.

{\bf Subcase 4b.2.iv. $M_1, M_2 \gtrsim M_3^{\frac{3}{N-1}}$.} By H\"older's inequality and Lemma $\ref{lem-str-est-Y}$, we have
\begin{align*}
\Big| \int_0^T \int_{\R^N} v_1 z_2 \scal{\nabla}^{\frac{N-4}{2}} z_3 v_4 dx dt\Big| &\leq \|v_1\|_{L^{\frac{2(N+4)}{N}}_{t,x}} \|z_2\|_{L^{\frac{N+4}{2}}_{t,x}} \|\scal{\nabla}^{\frac{N-4}{2}} z_3\|_{L^{\frac{N+4}{2}}_{t,x}} \|v_4\|_{L^{\frac{2(N+4)}{N}}_{t,x}} \\
&\lesssim M_1^{-\frac{N-4}{2}} M_2^{-\gamma} M_3^{\frac{N-4}{2}-\gamma} \|v_1\|_{X^{\frac{N-4}{2}}} \prod_{j=2}^3 \|\scal{\nabla}^\gamma z_j\|_{L^{\frac{N+4}{2}}_{t,x}} \|v_4\|_{Y^0}\\
&\lesssim M_3^{\frac{(N-4)^2}{2(N-1)} -\frac{(N+2)\gamma}{N-1}+}  \|v_1\|_{X^{\frac{N-4}{2}}} \prod_{j=2}^3 \|\scal{\nabla}^\gamma z_j\|_{L^{\frac{N+4}{2}}_{t,x}} \\
&\lesssim T^{0+} R^2 \|v_1\|_{X^{\frac{N-4}{2}}}
\end{align*}
outside a set of probability at most
\[
C \exp\left(-cR^2T^{-\frac{4}{N+4}+} \|f\|^{-2}_{H^\gamma}\right)
\]
as long as  
\[
\gamma >\frac{(N-4)^2}{2(N+2)}
\]
which is again satisfied if \eqref{cond-gamma-proof} holds.

Collecting the above cases, we prove \eqref{proba-non-est-1}. 

$\bullet$ We next estimate \eqref{global-proba-non-est-1}, the estimate \eqref{global-proba-non-est-2} is treated in a similar manner. Given $M\geq 1$, we define
\[
\tilde{\Phi}_M(v(t)):= \mp i \int_0^t U_\mu(t-s) P_{\leq M} \Nc(v+\vareps z)(s) ds.
\]
By Bernstein's inequality and H\"older's inequality, we have
\begin{align*}
\|P_{\leq M} \Nc(v+\vareps z)\|_{L^1_t H^{\frac{N-4}{2}}_x(\R \times \R^N)} &\lesssim M^{\frac{N-4}{2}} \|\Nc(v+\vareps z)\|_{L^1_t L^2_x(\R \times \R^N)} \\
&\lesssim M^{\frac{N-4}{2}} \|v\|^3_{L^3_t L^6_x(\R \times \R^N)} + M^{\frac{N-4}{2}} \vareps \|z\|^3_{L^3_tL^6_x(\R \times \R^N)}.
\end{align*}
Thanks to the global in time probabilistic Strichartz estimates and the fact $\left(3, \frac{6N}{3N-8}\right)$ is Biharmonic admissible with $6\geq \frac{6N}{3N-8}$, we see that the second term in the right hand side is finite almost surely. On the other hand, by Sobolev embedding and Lemma $\ref{lem-str-est-Y}$, we have
\[
\|v\|_{L^3_tL^6_x(\R \times \R^N)} \lesssim \||\nabla|^{\frac{N-4}{3}} v\|_{L^3_t L^{\frac{6N}{3N-8}}_x(\R \times \R^N)} \lesssim \|v\|_{X^{\frac{N-4}{2}}(\R)} <\infty.
\]
This shows that for each $M\geq 1$, $P_{\leq M} \Nc(v+\vareps z) \in L^1_t H^{\frac{N-4}{2}}(\R\times \R^N)$ almost surely. Thus, by Lemma $\ref{lem-lin-est}$,
\begin{align} \label{est-tilde-Phi-M}
\|\tilde{\Phi}_M(v)\|_{X^{\frac{N-4}{2}}(\R)} \lesssim \sup_{v_4 \in Y^0(\R) \atop \|v_4\|_{Y^0}=1} \left|\int_\R \int_{\R^N} \scal{\nabla}^{\frac{N-4}{2}}[\Nc(v+\vareps z)(t,x)] \overline{v}_4(t,x) dx dt \right|
\end{align}
almost surely, where $v_4 = P_{\leq M} v_4$. As above, we will estimate the right hand side of \eqref{est-tilde-Phi-M} independent of the cutoff size $M\geq 1$, by performing a case-by-case analysis of expressions of the form
\begin{align} \label{est-tilde-Phi-M-w}
\left|\int_{\R} \int_{\R^N} \scal{\nabla}^{\frac{N-4}{2}}(w_1w_2w_3)v_4 dx dt\right|,
\end{align}
where $\|v_4\|_{Y^0(\R)} \leq 1$ and $w_j=v$ or $z, j=1,2,3$. Then, letting $M\rightarrow \infty$, the same estimate holds for $\tilde{\Phi}(v)$ without any cutoff. 

The rest of the proof follows in a similar manner as the proof of the first part by changing the time interval from $[0,T)$ to $\R$ and replacing $z$ by $\vareps z$. Note that $\left(\frac{N+4}{2}, \frac{2N(N+4)}{N^2+4N-16}\right)$, $\left(4, \frac{2N}{N-2}\right)$ and $\left(\frac{6(N+4)}{N+8},\frac{6N(N+4)}{3N^2+8N-32}\right)$ are Biharmonic admissible and $\frac{N+4}{2} \geq \frac{2N(N+4)}{N^2+4N-16}, 4\geq \frac{2N}{N-2}$ and $\frac{6(N+4)}{N+8} \geq \frac{6N(N+4)}{3N^2+8N-32}$, we can use the global in time probabilistic Strichartz estimates for $(q,\tilde{r})=\left(\frac{N+4}{2},\frac{N+4}{2}\right)$, $(4,4)$ and $\left(\frac{6(N+4)}{N+8},\frac{6(N+4)}{N+8}\right)$, for instance
\[
\|\scal{\nabla}^\gamma z\|_{L^{\frac{N+4}{2}}_{t,x}} \leq \frac{R}{\vareps}
\]
outside a set of probability at most
\[
C \exp\left(-cR^2\vareps^{-2}\|f\|^{-2}_{H^\gamma}\right).
\]
We see that the contribution to \eqref{est-tilde-Phi-M-w} is given by 
\[
\text{\bf Case 2:} \quad R^3, \quad \text{\bf Case 3:} \quad R \prod_{j=1}^2 \|v_j\|_{X^{\frac{N-4}{2}}(\R)}, \quad \text{\bf Case 4:} \quad R^2 \|v_1\|_{X^{\frac{N-4}{2}}(\R)}
\]
outside a set of probability at most
\[
C \exp\left(-cR^2\vareps^{-2}\|f\|^{-2}_{H^\gamma}\right)
\]
in all cases as long as $\gamma>\gamma_N$. The proof is complete.
\end{proof}

\section{Probabilistic well-posedness}
\label{S4}
\setcounter{equation}{0}

We are now able to prove the almost sure local well-posedness for \eqref{4NLS} given in Theorem $\ref{theo-almo-sure-lwp}$.

\noindent {\bf Proof of Theorem $\ref{theo-almo-sure-lwp}$.} We will show that \eqref{4NLS-v} is almost sure locally well-posed. To this end, we consider
\[
\mathcal{X}:= \left\{ v \in X^{\frac{N-4}{2}}([0,T]) \cap C([0,T], H^{\frac{N-4}{2}}(\R^N)) \ : \ \|v\|_{X^{\frac{N-4}{2}}([0,T])} \leq \eta\right\}
\]
equipped with the distance
\[
d(v_1,v_2):= \|v_1-v_2\|_{X^{\frac{N-4}{2}}([0,T])}
\]
for some $T, \eta>0$ to be chosen later. It is enough to show that the functional 
\[
\Phi(v(t))= \mp i \int_0^t U_\mu(t-s) \Nc(v+z)(s) ds
\]
is a contraction on $(\mathcal{X},d)$. By Proposition $\ref{prop-proba-non-est}$, we have for $0<T\leq 1$, there exists $0<\vartheta \ll 1$ such that
\begin{align*}
\|\Phi(v)\|_{X^{\frac{N-4}{2}}([0,T])} &\leq C_1\left( \|v\|^3_{X^{\frac{N-4}{2}}([0,T]} + T^\vartheta R^3\right), \\
\|\Phi(v_1)-\Phi(v_2)\|_{X^{\frac{N-4}{2}}([0,T])} &\leq C_2 \left( \sum_{j=1}^2 \|v_j\|^2_{X^{\frac{N-4}{2}}([0,T])} + T^\vartheta R^2 \right) \|v_1-v_2\|_{X^{\frac{N-4}{2}}([0,T])}
\end{align*}
for all $v,v_1,v_2 \in X^{\frac{N-4}{2}}([0,T])$ and $R>0$, outside a set of probability at most
\[
C \exp\left(-cR^2\|f\|^{-2}_{H^\gamma(\R^N)}\right).
\]
It follows that for $v,v_1,v_2 \in \mathcal{X}$, 
\begin{align*}
\|\Phi(v)\|_{X^{\frac{N-4}{2}}([0,T])} &\leq C_1 \left(\eta^3+ T^\vartheta R^3 \right), \\
d(\Phi(v_1),\Phi(v_2)) &\leq C_2 \left(2\eta^2+T^\vartheta R^2\right) d(v_1,v_2)
\end{align*}
outside a set of probability at most
\[
C \exp\left(-cR^2\|f\|^{-2}_{H^\gamma(\R^N)}\right).
\]
We choose $\eta>0$ small so that
\[
C_1 \eta^2 \leq \frac{1}{2}, \quad 2C_2 \eta^2 \leq\frac{1}{4}.
\]
For given $R\gg 1$, we choose $T=T(R)$ so that
\[
C_1T^\vartheta R^3 \leq \frac{\eta}{2}, \quad C_2 T^\vartheta R^2 \leq \frac{1}{4},
\]
hence
\[
T^\vartheta = \min \left\{\frac{\eta}{2C_1R^3}, \frac{1}{4C_2 R^2} \right\}.
\]
With such choices, we see that $\Phi$ is a contraction on $(\mathcal{X},d)$ outside a set of probability at most
\[
C \exp \left(-cR^2\|f\|^{-2}_{H^\gamma(\R^N)}\right) \sim C \exp \left(-cT^{-\theta} \|f\|^{-2}_{H^\gamma(\R^N)}\right)
\]
for some $\theta >0$. The proof is complete.
\hfill $\Box$

We next prove the probabilistic small data global well-posedness and scattering for \eqref{4NLS} given in Theorem $\ref{theo-proba-small-gwp}$. 

\noindent {\bf Proof of Theorem $\ref{theo-proba-small-gwp}$.}
We consider
\[
\mathcal{Y}: = \left\{ v\in X^{\frac{N-4}{2}}(\R) \cap C(\R, H^{\frac{N-4}{2}}(\R^N)) \ : \ \|v\|_{X^{\frac{N-4}{2}}(\R)} \leq \delta \right\}
\]
equipped with the distance
\[
d(v_1,v_2):= \|v_1-v_2\|_{X^{\frac{N-4}{2}}(\R)}
\]
for some $\delta>0$ to be chosen later. For $0<\vareps \ll 1$, we will show that the functional
\[
\tilde{\Phi}(v(t)):= \mp i \int_0^t U_\mu(t-s) \Nc(v+ \vareps z)(s) ds
\]
is a contraction on $(\mathcal{Y},d)$. By Proposition $\ref{prop-proba-non-est}$ with $R=\delta$, we have for any $v,v_1,v_2 \in \mathcal{Y}$, 
\begin{align*}
\|\tilde{\Phi}(v)\|_{X^{\frac{N-4}{2}}(\R)} &\leq 2C_3 \delta^3, \\
d(\tilde{\Phi}(v_1),\tilde{\Phi}(v_2)) &\leq 3C_4 \delta^2 d(v_1,v_2)
\end{align*}
outside a set of probability at most
\[
C \exp \left(-c\delta^2\vareps^{-2}\|f\|^{-2}_{H^\gamma(\R^N)}\right).
\]
By choosing $\delta>0$ small so that
\begin{align} \label{defi-delta}
2C_3 \delta^2 \leq 1, \quad 3C_4\delta^2\leq \frac{1}{2},
\end{align}
we see that $\tilde{\Phi}$ is a contraction on $(\mathcal{Y},d)$ outside a set of probability at most 
\[
C\exp \left(-c\delta^2\vareps^{-2}\|f\|^{-2}_{H^\gamma(\R^N)}\right).
\]
Noting that $\delta$ is an absolute constant, we conclude that for each $0<\vareps\ll 1$, there exists a set $\Omega_\vareps \subset \Omega$ such that
\begin{itemize}
	\item $\Pc(\Omega \backslash \Omega_\vareps) \leq C \exp \left(-c\vareps^{-2} \|f\|^{-2}_{H^\gamma(\R^N)}\right)$;
	\item For each $\omega \in \Omega_\vareps$, there exists a unique global in time solution to \eqref{4NLS} with initial data $\vareps f^\omega$ in the class
	\[
	\vareps U_\mu(t) f^\omega + C(\R, H^{\gamc}(\R^N)) \subset C(\R, H^\gamma(\R^N)).
	\]
\end{itemize}
It remains to show the scattering. Fix $\omega \in \Omega_\vareps$ and let $v=v(\vareps, \omega)$ be the global in time solution to \eqref{4NLS-v} constructed above. We will show that there exists $f_+^\omega \in H^{\frac{N-4}{2}}(\R^N)$ such that
\begin{align} \label{est-v}
U_\mu(-t) v(t) = \mp i\int_0^t U_\mu(-s) \Nc(v+\vareps z)(s) ds \rightarrow f^\omega_+
\end{align}
in $H^{\frac{N-4}{2}}(\R^N)$ as $t\rightarrow \infty$. Set $w(t)=U_\mu(-t) v(t)$. For $0<t_1 \leq t_2 <\infty$, we have
\begin{align*}
U_\mu(t_2)(w(t_2) - w(t_1)) &= \mp i \int_{t_1}^{t_2} U_\mu(t_2-s) \Nc(v+\vareps z)(s) ds \\
&=\mp i \int_0^{t_2} U_\mu(t_2-s) \chi_{[t_1,\infty)}(s) \Nc(v+\vareps z)(s) ds=: I(t_1,t_2).
\end{align*}
In the following, we view $I(t_1,t_2)$ as a function on $t_2$ and estimate its $X^{\frac{N-4}{2}}([0,\infty))$-norm. We now revisit the computation in the proof of Proposition $\ref{prop-proba-non-est}$. 

In Case 1, we proceed slightly differently. By Lemma $\ref{lem-lin-est}$ and H\"older's inequality, we have
\begin{align}
\|I(t_1,t_2)\|_{X^{\frac{N-4}{2}}([0,\infty))} &\lesssim \sup_{v_4 \in Y^0([0,\infty)) \atop \|v_4\|_{Y^0}=1} \left| \int_0^\infty \int_{\R^N} \chi_{[t_1,\infty)}(t) \scal{\nabla}^{\frac{N-4}{2}} v \overline{v} v v_4 dx dt \right| \nonumber \\
&\leq \|\scal{\nabla}^{\frac{N-4}{2}} v\|_{L^{\frac{2(N+4)}{N}}_{t,x}([t_1,\infty)\times \R^N)} \|v\|^2_{L^{\frac{N+4}{2}}_{t,x}([t_1,\infty)\times \R^N)}. \label{est-case-1}
\end{align}
By Lemma $\eqref{lem-str-est-Y}$, we have
\[
\|\scal{\nabla}^{\frac{N-4}{2}} v\|_{L^{\frac{2(N+4)}{N}}_{t,x}(\R \times \R^N)} + \|v\|_{L^{\frac{N+4}{2}}_{t,x}(\R \times \R^N)} \lesssim \|v\|_{X^{\frac{N-4}{2}}(\R)} \leq \delta.
\]
Then by the monotone convergence theorem, \eqref{est-case-1} tends to 0 as $t_1 \rightarrow \infty$. 

In Cases 2, 3 and 4, we had at least one factor $z$. We multiply the cutoff function $\chi_{[t_1,\infty)}$ only on the $\vareps z$-factors but not on the $v$-factors. Note that $\|v\|_{X^{\frac{N-4}{2}}(\R)} \leq \delta$. As in the proof of Proposition $\ref{prop-proba-non-est}$, we estimate at least a small portion of these $z$-factors in $\|\scal{\nabla}^\gamma \vareps z^\omega\|_{L^r_{t,x}([t_1,\infty)\times \R^N)}$ with $q=\frac{N+4}{2}$ or 4 or $\frac{6(N+4)}{N+8}$. Recall that $\|\scal{\nabla}^\gamma \vareps z^\omega\|_{L^r_{t,x}([t_1,\infty)\times \R^N)} \leq \delta$ for $\omega \in \Omega_\vareps$. Hence, by the monotone convergence theorem, we have $\|\scal{\nabla}^\gamma \vareps z^\omega\|_{L^r_{t,x}([t_1,\infty))} \rightarrow 0$ as $t_1 \rightarrow \infty$. Thus, the contributions from Cases 2, 3 and 4 tend to 0 as $t_1\rightarrow \infty$. Therefore,
\[
\lim_{t_1\rightarrow \infty} \|I(t_1,t_2)\|_{X^{\frac{N-4}{2}}([0,\infty))} =0.
\] 
In conclusion, we obtain
\begin{align*}
\lim_{t_1\rightarrow \infty} \sup_{t_2>t_1} \|w(t_2)-w(t_1)\|_{H^{\frac{N-4}{2}}(\R^N)} &= \lim_{t_1\rightarrow \infty} \sup_{t_2>t_1} \|U_\mu(t_2)(w(t_2)-w(t_1))\|_{H^{\frac{N-4}{2}}(\R^N)} \\
&= \lim_{t_1\rightarrow\infty} \|I(t_1,t_2)\|_{L^\infty_{t_2} H^{\frac{N-4}{2}}_x([0,\infty)\times \R^N)} \\
&\lesssim \lim_{t_1\rightarrow \infty} \|I(t_1,t_2)\|_{X^{\frac{N-4}{2}}([0,\infty))} =0.
\end{align*}
This proves \eqref{est-v} and the scattering of $u^\omega(t) = \vareps U_\mu(t)f^\omega + v^\omega(t)$, which completes the proof of Theorem $\ref{theo-proba-small-gwp}$.
\hfill $\Box$

Finally, we prove the global well-posedness and scattering with a large probability given in Theorem $\ref{theo-gwp-large}$. We follow the idea of \cite{BOP-TAMS}, that is to exploit the dilation symmetry \eqref{scaling} of the cubic 4NLS \eqref{4NLS}. Denote
\[
f_\lambda(x):= \lambda^2 f(\lambda x), \quad \lambda>0.
\]
We have
\begin{align} \label{scaling-f}
\|f_\lambda\|_{\dot{H}^\gamma(\R^N)} = \lambda^{\gamma-\frac{N-4}{2}} \|f\|_{\dot{H}^\gamma(\R^N)}.
\end{align}
If $\gamma <\gamc=\frac{N-4}{2}$, 
then we can make the $H^\gamma$-norm of the scaled function $f_\lambda$ small by taking $\lambda \gg 1$. The issue is that the Strichartz estimates we employ in proving probabilistic well-posedness are (sub)critical and do not become small enven if we take $\lambda \gg 1$. It is for this reason that we consider the randomization $f^\omega_\lambda$ on dilated cubes. 

\noindent {\bf Proof of Theorem $\ref{theo-gwp-large}$.}
Fix $f \in H^\gamma(\R^N)$ with $\gamma \in (\gamma_N, \gamc)$, where $\gamma_N$ is as in \eqref{defi-gam-N}. Let $f^{\omega,\lambda}$ be its randomization on dilated cubes of scale $\lambda$ as in \eqref{defi-rando-lambda}. Instead of considering \eqref{4NLS} with $u_0 = f^{\omega,\lambda}$, we consider the scaled Cauchy problem
\begin{align} \label{4NLS-scale}
\left\{
\begin{array}{rcl}
i\partial_t u_\lambda - \Delta^2 u_\lambda &=& \pm |u_\lambda|^2 u_\lambda,  \\
u_\lambda|_{t=0} &=& u_{0,\lambda} = (f^{\omega,\lambda})_\lambda,
\end{array}
\right.
\end{align}
where $u_\lambda$ is as in \eqref{scaling} and $(f^{\omega,\lambda})_\lambda(x) = \lambda^2 f^{\omega,\lambda}(\lambda x)$ is the scaled randomization. For simplicity, we denote $(f^{\omega,\lambda})_\lambda$ by $f^{\omega,\lambda}_\lambda$ in the following. We denote the linear and nonlinear part of $u_\lambda$ by $z_\lambda(t) = z^\omega_\lambda(t) := U_0(t) f^{\omega,\lambda}_\lambda$ and $v_\lambda(t) := u_\lambda(t) - U_0(t) f^{\omega,\lambda}_\lambda$. We reduce \eqref{4NLS-scale} to 
\begin{align} \label{4NLS-scale-v}
\left\{
\begin{array}{rcl}
i\partial_t v_\lambda - \Delta^2 v_\lambda &=& \pm |v_\lambda + z_\lambda|^2 (v_\lambda+z_\lambda),   \\
v_\lambda|_{t=0} &=& 0.
\end{array}
\right.
\end{align}
Note that if $u$ satisfies \eqref{4NLS} with initial data $u(0)= f^{\omega,\lambda}$, then $u_\lambda, z_\lambda$ and $v_\lambda$ are the scalings of $u, z:=U_0(t) f^{\omega,\lambda}$ and $v:= u-z$ respectively. In fact, it is clear for $u_\lambda$ by using \eqref{scaling}. For $z_\lambda$ and $v_\lambda$, this follows from the following observation:
\begin{align} \label{scaling-lambda}
\mathcal{F}_x \left[ \left(U_0(t) f^{\omega,\lambda}\right)_\lambda \right](\xi) = \lambda^{2-N} e^{-i\lambda^4 t |\lambda^{-1} \xi|^4} \widehat{f^{\omega,\lambda}}(\lambda^{-1} \xi) = e^{-it|\xi|^4} \widehat{f^{\omega,\lambda}_\lambda} (\xi) = \widehat{z_\lambda}(t,\xi).
\end{align}
Define
\[
\tilde{\Phi}_\lambda(v_\lambda(t)) := \mp i \int_0^t U_0(t-s) \Nc(v_\lambda + z_\lambda)(s) ds.
\]
We will show that there exists $\lambda_0=\lambda_0(\vareps, \|f\|_{H^\gamma(\R^N)})>0$ such that, for $\lambda >\lambda_0$, the estimates \eqref{global-proba-non-est-1} and \eqref{global-proba-non-est-2} in Proposition $\ref{prop-proba-non-est}$ (with $\tilde{\Phi}$ replaced by $\tilde{\Phi}_\lambda$) hold with $R= \delta$ outside a set of probability strictly smaller than $\vareps$, where $\delta$ is as in \eqref{defi-delta}. In fact, we first observe that
\[
\psi(D-n) f_\lambda = \left( \psi_\lambda(D-\lambda^{-1} n) f\right)_\lambda.
\]
Hence, we have
\begin{align} \label{scaling-rando}
f^{\omega,\lambda}_\lambda = (f^{\omega,\lambda})_\lambda = \sum_{n\in \Z^N} g_n(\omega) \psi(D-n) f_\lambda.
\end{align}
Given $\delta$ as in \eqref{defi-delta} and $\lambda>0$, we define
\[
\Omega_{1,\lambda}:= \left\{ \omega \in \Omega \ : \ \|U_0(t) f^{\omega,\lambda}_\lambda\|_{L^r_t W^{\gamma,r}_x(\R \times \R^N)} \leq \delta, \quad r= \frac{N+4}{2}, 4, \frac{6(N+4)}{N+8} \right\}.
\]
We also define
\[
\Omega_{2,\lambda}:= \left\{ \omega \in \Omega \ : \ \|f^{\omega,\lambda}_\lambda\|_{H^\gamma(\R^N)} \leq \delta \right\}.
\]
Now, set $\Omega_\lambda= \Omega_{1,\lambda} \cap \Omega_{2,\lambda}$. It follows from \eqref{scaling-rando}, Lemma $\ref{lem-differ}$ and Lemma $\ref{lem-glo-proba-str-est}$ that
\[
\Pc(\Omega \backslash \Omega_\lambda) \leq C \exp \left(-c\delta^2\|f_\lambda\|^{-2}_{H^\gamma(\R^N)}\right) \leq C \exp \left(-c \delta^2\lambda^{-2\gamma+N-4} \|f\|^{-2}_{H^\gamma(\R^N)}\right)
\]
for $\lambda \geq 1$. Note that
\begin{align*}
\|f_\lambda\|^2_{H^\gamma} &= \|f_\lambda\|^2_{L^2} + \|f_\lambda\|^2_{\dot{H}^\gamma} \\
&= \lambda^{-N+4} \|f\|^2_{L^2} + \lambda^{2\gamma-N+4} \|f\|^2_{\dot{H}^\gamma} \\
&= \lambda^{2\gamma-N+4} \left( \lambda^{-2\gamma} \|f\|^2_{L^2} + \|f\|^2_{\dot{H}^\gamma}\right) \\
&\leq \lambda^{2\gamma-N+4} \|f\|^2_{H^\gamma}
\end{align*}
since $\lambda \geq 1$. By setting
\[
\lambda_0 \sim \left( \frac{\log \left(\frac{1}{\vareps}\right) \|f\|^2_{H^\gamma(\R^N)}}{\delta^2}\right)^{\frac{1}{N-4-2\gamma}},
\]
we have
\[
\Pc(\Omega \backslash \Omega_\lambda) <\vareps
\]
for all $\lambda >\lambda_0$. Note that $\lambda_0 \rightarrow \infty$ as $\vareps \rightarrow 0$. 	

Recall that the pairs $\left(\frac{N+4}{2}, \frac{N+4}{2}\right), (4,4)$ and $\left(\frac{6(N+4)}{N+8}, \frac{6(N+4)}{N+8}\right)$ are the only relevant values of the space-time Lebesgue indices controlling the random forcing term in the proof of Proposition $\ref{prop-proba-non-est}$. Hence, the estimates \eqref{global-proba-non-est-1} and \eqref{global-proba-non-est-2} in Proposition $\ref{prop-proba-non-est}$ (with $\tilde{\Phi}$ replaced by $\tilde{\Phi}_\lambda$) hold with $R=\delta$ for each $\omega \in \Omega_\lambda$. Then repeating the proof of Theorem $\ref{theo-proba-small-gwp}$, we see that for each $\omega \in \Omega_\lambda$, there exists a unique global solution $u_\lambda$ to \eqref{4NLS-scale} which scatters both forward and backward in time. By undoing the scaling, we obtain a unique global solution $u$ to \eqref{4NLS} with initial data $f^{\omega,\lambda}$ for each $\omega \in \Omega_\lambda$. Moreover, scattering for $u_\lambda$ implies the scattering for $u$. Indeed, as in Theorem $\ref{theo-proba-small-gwp}$, there exists $f^\omega_{+,\lambda} \in H^{\frac{N-4}{2}}(\R^N)$ such that
\[
\lim_{t\rightarrow \infty} \|u_\lambda(t) - U_0(t) f^{\omega,\lambda}_\lambda - U_0(t) f^\omega_{+,\lambda} \|_{H^{\frac{N-4}{2}}(\R^N)} =0.
\]
A computation similar to \eqref{scaling-lambda} gives
\[
U_0(t) f^{\omega,\lambda}_\lambda + U_0(t) f^\omega_{+,\lambda} = \left(U_0(t) f^{\omega,\lambda} + U_0(t) f^\omega_+ \right)_\lambda,
\]
where $f^\omega_+ = (f^\omega_{+,\lambda})_{\lambda^{-1}} \in H^{\frac{N-4}{2}}(\R^N)$. Then, by \eqref{scaling-f}, we obtain
\[
\lim_{t\rightarrow \infty} \|u(t) - U_0(t) f^{\omega,\lambda} - U_0(t) f^\omega_+ \|_{H^{\frac{N-4}{2}}(\R^N)} =0.
\]
This proves that $u$ scatters forward in time. The proof is complete.
\hfill $\Box$

	\section*{Acknowledgement}
	This work was supported in part by the Labex CEMPI (ANR-11-LABX-0007-01). The author would like to express his deep gratitude to his wife - Uyen Cong for her encouragement and support. He also would like to thank the reviewers for their corrections and valuable suggestions.

\end{document}